\input ./style/arxiv-general.cfg
\documentclass[MSNbibl,number,citesort,seceqn,dvips]{arxbj}
\makeatletter
   \@ifpackageloaded{graphicx}{}{\usepackage{graphicx}}
\makeatother
\usepackage{mathbh}

\volume{22}
\issue{2}
\pubyear{2016}
\firstpage{681}
\lastpage{710}
\doi{10.3150/14-BEJ672} 
\docsubty{FLA}

\makeatletter
\newtheorem{theorem}{Theorem}[section]
\newtheorem{proposition}[theorem]{Proposition}
\newtheorem{lemma}[theorem]{Lemma}
\newtheorem{lemmaa}{Lemma}[section]
\newproclaim{definition}[theorem]{Definition}
\newremark{remark}[theorem]{Remark}

\def\R{\mathbb{R}}
\def\eps{{\varepsilon}}
\def\un{\mathbh{1}}

\newcommand{\bx}{\mathbf{x}}
\newcommand{\by}{\mathbf{y}}
\newcommand{\bz}{\mathbf{z}}
\newcommand{\bv}{\mathbf{v}}

\newcommand{\medm}{{\delta'}}
\newcommand{\medp}{{\delta}}
\def\phi{{\varphi}}

\def\CC{\mathcal C}
\def\EE{\mathbf E}

\def\Var{\operatorname{Var}}

\newcommand{\Es}{{\mathsf{E}}}
\newcommand{\Ff}{\mathcal{F}}
\newcommand{\Ps}{\mathsf{P}}
\makeatother

\begin{document}
\begin{frontmatter}

\title{Long time behavior of stochastic hard ball~systems}
\runtitle{Long time behavior of stochastic hard ball systems}

\begin{aug}
\author[A]{\inits{P.}\fnms{Patrick}~\snm{Cattiaux}\thanksref{A}\ead[label=e1]{cattiaux@math.univ-toulouse.fr}},
\author[B]{\inits{M.}\fnms{Myriam}~\snm{Fradon}\thanksref{B}\ead[label=e2]{Myriam.Fradon@univ-lille1.fr}},
\author[C]{\inits{A.~M.}\fnms{Alexei M.}~\snm{Kulik}\thanksref{C}\ead[label=e3]{kulik@imath.kiev.ua}}\\ 
\and
\author[D]{\inits{S.}\fnms{Sylvie}~\snm{Roelly}\corref{}\thanksref{D}\ead[label=e4]{roelly@math.uni-potsdam.de}}
\address[A]{Institut de Math\'ematiques de Toulouse, CNRS UMR 5219,
Universit\'e de Toulouse, 118 route de Narbonne, 31062 Toulouse Cedex 09,
France. \printead{e1}}
\address[B]{U.F.R. de Math\'ematiques, CNRS UMR 8524, Universit\'e
Lille 1, 59655 Villeneuve d'Ascq Cedex, France. \printead{e2}}
\address[C]{Institute of
Mathematics, Ukrainian National Academy of Sciences, Kyiv 01601, Tereshchenkivska str. 3, Ukraine. \printead{e3}}
\address[D]{Institut f\"ur Mathematik der Universit\"at Potsdam, Am
Neuen Palais 10, 14469 Potsdam, Germany. \printead{e4}}
\end{aug}

\received{\smonth{2} \syear{2014}}
\revised{\smonth{7} \syear{2014}}


\begin{abstract}
We study the long time behavior of a system of $n=2,3$ Brownian hard
balls, living in $\R^d$ for $d \ge2$, submitted to a mutual attraction
and to elastic
collisions.
\end{abstract}

%
\end{frontmatter}

\section{Introduction and main results}\label{Intro}

Consider $n$ hard balls with radius $r/2$ and centers $X_1,\ldots,X_n$
located in $\R^d$ for some $d\ge2$. They are moving randomly and when they
meet, they are performing elastic collisions. We are interested in the
long time behavior of such a dynamics, where the centers of the balls
are moving
according to a Brownian motion in a Gaussian type pair potential. It is
modelized by the following system of
stochastic differential equations with
reflection
\[
(\mathrm{A}) \cases{ %
\mbox{for } i \in\{1,\ldots,n\}, t
\in\R^+ ,
\cr
\displaystyle X_i(t) = X_i(0) + W_i(t) - a \sum
_{j=1}^n \int_0^t
\bigl(X_i(s)-X_j(s)\bigr) \,\mathrm{d}s
\cr
\hspace*{33pt}{} \displaystyle + \sum
_{j=1}^n \int_0^t
\bigl(X_i(s)-X_j(s)\bigr) \,\mathrm{d}L_{ij}(s) ,
\cr
L_{ij}(0) = 0 , L_{ij} \equiv L_{ji}
\mbox{ and }\displaystyle L_{ij}(t) = \int_0^t
\un_{|X_i(s)-X_j(s)|=r} \,\mathrm{d}L_{ij}(s), L_{ii} \equiv0,}
\]
where $W_1,\ldots,W_n$ are $n$ independent standard Wiener processes.
The local time $L_{ij}$ describes the elastic collision (normal mutual
reflection) between
balls $i$ and $j$.
The parameter $a$ is assumed to be non-negative. Therefore, the drift
term derives from an attractive quadratic potential.

Note that the Markov process $X$ satisfying (A) admits a unique (up to
a multiplicative constant) unbounded invariant measure
$\mu_a$ defined on $(\mathbb R^d)^n$ by:
%
\begin{equation}
\label{eqmesure} \mathrm{d}\mu_a(\bx)= \mathrm{e}^{- a \sum_{i,j} |x_i-x_j|^2 /2} \un
_{D}(\bx)
\,\mathrm{d}\bx.
\end{equation}

Here $\bx= (x_1,\ldots,x_n) \in(\mathbb R^d)^n$ and $D$ is the interior
of the set of allowed configurations, that is,
%
\begin{equation}
\label{eqdefD} D = \bigl\{\bx\in\bigl(\mathbb R^d\bigr)^n;
|x_i - x_j|>r \mbox{ for all } i \neq j \bigr\}.
\end{equation}
Clearly the measure $\mu_a$ is invariant under the simultaneous
translations of the $n$ balls, that is under
any transformation of the form $(x_1,\ldots,x_n) \mapsto
(x_1+u,\ldots,x_n+u), u \in\mathbb R^d$.

Indeed we are even more interested by the intrinsic dynamics of the
system, that is, by the system of balls viewed from their center of
mass, called $G:=\frac{1}n
(X_1+ \cdots+X_n)$.
This (fictitious) point undergoes a Brownian motion in $\mathbb R^d$
with covariance $\frac{1}{n}$ Id (notice the absence of reflection term).
Choosing $G$ as the (moving) origin of the ambient space $\R^d$, we
therefore consider the process $Y$ of the relative positions,
$Y_i=X_i-G, i=1, \ldots,n$,
which satisfies
\[
(\mathrm{B}) \cases{ %
\mbox{for } i \in\{1,\ldots,n\} , t
\in\R^+ ,
\cr
\displaystyle Y_i(t) = Y_i(0) + M_i(t) - a \sum
_{j=1}^n \int_0^t
\bigl(Y_i(s)-Y_j(s)\bigr) \,\mathrm{d}s + \sum
_{j=1}^n \int_0^t
\bigl(Y_i(s)-Y_j(s)\bigr) \,\mathrm{d}L_{ij}(s),
\vspace*{2pt}\cr
L_{ij}(0) = 0 , L_{ij} \equiv L_{ji}
\mbox{ and } \displaystyle L_{ij}(t) = \int_0^t
\un_{|Y_i(s)-Y_j(s)|=r} \,\mathrm{d}L_{ij}(s), L_{ii} \equiv0,}
\]
where the martingale term $(M_1,\ldots, M_n)$ is a new Brownian motion
with covariation
$
\langle M_i, M_k \rangle(t)=
(\frac{n-1}{n} \delta_{\{i=k\}} - \frac{1}n \delta_{\{
i\neq k\}} ) t \mathrm{Id} $.

The $(\R^{d})^n$-valued Markov process $Y(t)$ admits as unique
invariant probability measure
%
\begin{equation}
\label{eqmesurerestreinte} \mathrm{d}\pi_a(\by)= Z_a^{-1}
\mathrm{e}^{- a \sum_{i,j} |y_i-y_j|^2/2} \un_{ D'}(\by) \,\mathrm{d}\by
\end{equation}
for a well-chosen normalization constant $ Z_a$.
The domain $D'$, support of $\pi_a$ obtained as linear transformation
of $D$, is the following unbounded set
%
\begin{equation}
\label{eqdefD'} D':= \Biggl\{\by\in\bigl(\mathbb R^d
\bigr)^n ; |y_i - y_j|>r \mbox{ for all } i
\neq j , \sum_{i=1}^n y_i =0
\Biggr\}.
\end{equation}

Our aim is to describe the long time behavior of the process $Y$, that
is of the system of balls viewed from their center of mass.

Before explaining in more details the contents of the paper, let us
give an account of the existing literature and of related problems.

Existence and uniqueness of a strong solution for system (A) was first
obtained in Saisho and Tanaka \cite{ST86} with $a=0$. Extensions to $a
>0$ and to
$n=+\infty$ are done in
Fradon and R\oe lly \cite{FR1,FR2}, Fradon,
R\oe lly and Tanemura \cite{FRT}. Random radii $r$ were also studied in
Fradon \cite
{FradonGlob}, Fradon and Roelly \cite{FradonRoellyGlob}. The invariant
(in fact
reversible) measure
for the system is
discussed in Saisho and Tanaka \cite{ST87} and Fradon and R\oe lly \cite
{FR3} for an infinite number of balls.

The construction of the stationary process (i.e., starting from the
invariant measure) can also be performed by using Dirichlet forms
theory. Actually, $D$
intersected with any ball $B(0,R)\subset(\mathbb R^d)^n$ is a
Lipschitz domain (see the \hyperref[app]{Appendix}) so that one can
use results in
Bass and Hsu \cite{BassHsu}, Chen, Fitzsimmons and Williams \cite
{CFW}, Fukushima and Tomisaki \cite{Fuktom} to
build the Hunt process naturally associated to the Dirichlet form (see,
e.g., Fukushima, Oshima and Takeda \cite{FOT} for the
theory of Dirichlet forms)
%
\begin{equation}
\label{eqdir1} \mathcal E^R_a(f) = \int
_{D\cap B(0,R)} |\nabla f|^2 \,\mathrm{d}\mu_a.
\end{equation}
It is then enough to let $R$ go to infinity and show conservativeness
of the obtained process which is equivalent to non-explosion. This is standard.

The solution of (A) built by using stochastic calculus do coincide with
the Hunt process associated to the Dirichlet form $\mathcal E_a$
obtained for
$R=+\infty$.
Some properties, like the decomposition of the boundary into a
non-polar and a polar parts or the Girsanov's like structure are
discussed in Chen \textit{et al.} \cite{CFTYZ}.

For an infinite number of balls, such a construction is performed in
Osada \cite{Osada96}, Tanemura \cite{Tanemura96,Tanemura97}.

Let us recall also some regularity of the processes and their
associated semi-groups which we will need in the sequel.
For $x \in\bar D$, we denote by $P_t(x,\mathrm{d}y)$ the transition kernel of
the process $X(\cdot)$ starting from $x$ at time $t$.
It is well known that for all $x \in D$, $P_t(x,\mathrm{d}y)$ is absolutely
continuous with respect to the Lebesgue measure restricted to $\bar D$
(in particular does not charge the boundary $\partial D$). In addition
the density $p_t(x,y)$ is smooth as a function of the two variables
$x$ and $y$ in $D \times D$. This follows from standard elliptic
estimates (as in Bass and Hsu \cite{BassHsu2}) or from the use of Malliavin
calculus as
explained in Cattiaux \cite{Cx86,Cx92}. Furthermore, this density
kernel extends
smoothly up to the smooth part of the boundary (see Cattiaux \cite{Cx87,Cx92}).
But since the domain is (locally) Lipschitz, the potential theoretic
tools of Bass and Hsu \cite{BassHsu2} Sections 3 and 4 can be used to
show that
$(t,x,y) \mapsto p_t(x,y)$ extends continuously to $\mathbb R^+\times
\bar D \times\bar D$.
Actually Section~4 in Bass and Hsu \cite{BassHsu2} is written for
bounded Lipschitz
domain but extends easily to our situation by
localizing the Dirichlet form as we mentioned earlier and using
conservativeness (of course the function is no more uniformly continuous).
In particular, the process is Feller (actually strong Feller thanks to
Fukushima and Tomisaki \cite{Fuktom}).

Comparison with the killed process at the boundary shows that for any
$t>0$ and any starting $x \in D$, $p_t(x,y)>0$ for any $ y \in D$
(see, e.g., (3.15) and (3.16) in Bass and Hsu \cite{BassHsu2} and use
repeatedly the
Chapman--Kolmogorov relation to extend the result to all $t$ and $y$
introducing a chaining from $x$ to $y$). The previous continuity thus
implies that for all $t>0$, all compact subsets $K$ and $K'$ of $\bar D$
there exists a constant $C(t,K,K')>0$ such that
%
\begin{equation}
\label{positive} \mbox{for all $x\in K$ and $y\in K'$}\qquad
p_t(x,y) \geq C\bigl(t,K,K'\bigr).
\end{equation}
In particular compact sets are ``petite sets'' in the Meyn--Tweedie
terminology (Meyn and Tweedie~\cite{MT}) and for any compact set $K
\subset\bar D$ and
any $t>0$, the
%
\begin{equation}
\label{doeblin} (\textit{Local Dobrushin condition}) \qquad\sup_{x,x' \in
K}
\bigl\Vert P_t(x,\mathrm{d}y) - P_t\bigl(x',\mathrm{d}y\bigr)\bigr\Vert
_{\mathrm{TV}} < 2 ,
\end{equation}
is fulfilled, where $\Vert\cdot\Vert_{\mathrm{TV}}$ denotes the total
variation distance.

Another classical consequence is the uniqueness of the invariant
measure (since all invariant measures are actually equivalent) up to a
multiplicative constant.

Since $Y$ is deduced from $X$ by a smooth linear transformation,
similar statements are available for $Y$ in $D'$.
In particular, $Y$ is a Feller process satisfying the local Dobrushin
condition with a unique invariant probability measure.

Looking at long time behavior of such systems is not only interesting
by itself but relates, as $a \to+\infty$ (low temperature regime in
statistical
mechanics), to the following \emph{finite packing problem}: what is the
shape of a cluster of $n$ spheres -- with equal radii $r/2$ -- minimizing
their quadratic
energy,
that is, their second moment about their center of mass. (For a review of
different questions on finite packing, see the recent monograph B\"or\"
oczky \cite{Boeroe}.)
This problem, in spite of its simple statement and its numerous useful
applications, remains mainly open. Even for $d=2$ (so called \emph
{penny-packings}), only
the case $n \le7$ was solved by Temesv\'ari~\cite{Tem}.
For more
pennies, the optimal configurations are known only among the specific
class of hexagonal
packings (Chow \cite{Chow}). For $d=3$, one finds in Sloane \textit{et
al.} \cite{SHDC} a description
of the putatively optimal arrangements until $n \le32$.
For the case of infinitely many spheres and their celebrated densest
packing, we refer to Conway and Sloane \cite{ConwaySloane} or to
Fradon and R\oe lly \cite{FR3}, pages~99--100,
for recent references
and a more complete discussion.

Indeed, as $a \to+\infty$, the invariant measure $\pi_a$ concentrates
on the set of configurations with minimal quadratic energy that is, the set
\[
\EE_{\mathrm{min}} = \biggl\{\by\in D' ; V(\by):=\sum
_{i,j} |y_i-y_j|^2 =
\inf_{\bz\in\mathcal D'} \sum_{i,j}
|z_i-z_j|^2\biggr\} ,
\]
which obviously depends on $n$, $r$ and $d$. So looking simultaneously
at large $t$ and large $a$ furnishes some \emph{simulated annealing
algorithm} for the
uniform measure on $\EE_{\mathrm{min}}$ (see Theorem~\ref{th_exp_erg_3} for the
case $n=3$).

A similar (but different) algorithmic point of view is discussed in the
recent paper Diaconis, Lebeau and Michel \cite{DLM}, where the problem
under discussion is:
\textit{how can we place \emph{randomly } $n$ hard balls of radius
$r$ in
a given large ball (or hypercube)?} According to the introduction of
Diaconis, Lebeau and Michel \cite{DLM} this
problem is the origin of Metropolis algorithm. The authors relate the
asymptotics of the spectral gap of a discrete Metropolis algorithm
to the first Neumann eigenvalue (called $\nu_1$) for the Laplace
operator in $D'$ intersected with a large hypercube (see their Theorem~4.6).

There are several methods to attack the study of long time behavior for
Markov processes. In this paper, we will restrict ourselves to
exponential (or
geometric) ergodicity. Moreover, we will try to give some controls on
the rate of exponential ergodicity. Let us first recall some definitions.

\begin{definition}\label{defergod}
A Markov process $Z$ with transition distribution $P_t$ and invariant
measure $\pi$ is said to be exponentially ergodic if there exists
$\beta>0$ such that for
all initial condition $z$,
\[
\bigl\Vert P_t(z,\cdot) - \pi\bigr\Vert_{\mathrm{TV}} \leq C(z)
\mathrm{e}^{-\beta
t}.
\]
\end{definition}

If the function $z \mapsto C(z)$ is $\mu$-integrable, the previous
extends to any initial distribution $\mu$.

Our main result reads as follows.

\begin{theorem}\label{thmmain}
Consider a system of $n$ hard balls in $\mathbb R^d$ submitted to the
dynamics described by \emph{(A)}.
If $n=2,3$, the process $Y$ of their relative positions viewed from
their center of mass, described by the system \emph{(B)}, is
exponentially ergodic.
\end{theorem}

Remark that, if we are only interested in the convergence of the ball
system to the set of configurations with minimal energy in the large attraction
regime, the quantity of interest reduces to the $(\mathbb
R^+)^n$-valued system of the distances between the
centers of the $n$ balls and their center of mass. Its rate of
convergence to equilibrium is much faster that those of the $(\mathbb
R^d)^n$-valued process
$Y$.
For two balls, the difference is explicit when comparing Theorems \ref
{thm2} and \ref{thm2bis} in the next section.

Exponential ergodicity is connected to the existence of an exponential
coupling, as explained in Kulik \cite{Kul1,Kul2}, and is strongly dependent
on the existence of
exponential moments for the hitting time of compact subsets.
This method can be traced back to Veretennikov \cite{Veret}.
Let us give a precise statement taken from Kulik \cite{Kul1}, Theorem~2.2.

\begin{theorem}\label{thmalex}
Suppose that the process $Z$ satisfies the local Dobrushin condition.
Assume that we can find a real valued function $\Phi$ and a compact set
$K $ and positive
constants $c,\alpha$ such that:
\begin{longlist}[1.]
\item[1.]$\Phi$ is larger than 1 and $\Phi(z) \to+\infty$ as $|z|\to
+\infty$,
\item[2.] there exists $\alpha>0$ and $c>0$ such that for all initial
condition $z$,
\[
\Es_z\bigl(\Phi\bigl(Z(t)\bigr) \un_{\tau_K>t}\bigr) \leq c
\mathrm{e}^{-\alpha t} \Phi(z),
\]
where $\tau_K$ denotes the hitting time of $K$,
\item[3.]$\sup_{z\in K , t>0} \Es_z(\Phi(Z(t)) \un_{\Phi
(Z(t))>M}) \to0$ as $M \to+\infty$.
\end{longlist}
Then the process $Z$ is exponentially ergodic.
\end{theorem}

Though the result is only stated in the case $c=1$ in Kulik \cite
{Kul1}, the
method extends to any $c>0$ without difficulty.
It is instructive to compare this result with other forms of \emph
{Harris-type} theorems for exponential ergodicity, frequently used in
the literature.
Typically such theorems assume an irreducibility of the process on a
set $K$ (in our case this is the local Dobrushin condition) and its recurrence.
Recurrence assumption is formulated usually in the terms of the
generator $\mathcal{L}$ of the process, like the Foster--Lyapunov
condition, see, for example, Meyn and Tweedie \cite{MT},
%
\begin{equation}
\label{Foster} \mathcal L\Phi\leq-\alpha\Phi+ C \un_K.
\end{equation}
In Theorem~\ref{thmalex}, assumptions (2) and (3) can be interpreted as
a recurrence assumption, but since the generator
$\mathcal{L}$
is not involved therein, we call it an \emph{integral} Lyapunov
condition, while~(\ref{Foster}) is a \emph{differential} one.
In our framework, because of the presence of several local time terms
in (B), it is very difficult to find a Lyapunov function which
satisfy (\ref{Foster}) but we succeeded in showing that the quadratic
energy of the system, $V$, satisfies the more
tractable integral Lyapunov condition presented in Theorem~\ref{thmalex}.

It should be noted that, unfortunately, the exponential rate of
convergence $\beta$ in Theorem~\ref{thmalex} is difficult to express
explicitly in a compact
form, as it depends on $\alpha$ but also on other constants connected
with the behavior of the process,
in particular a quantitative version of the local Dobrushin condition
in $K$.
For an example of such an explicit expression, we refer the reader to
the end of
Section~3.2 in Kulik \cite{Kul2}. Similar formulae should appear in the
framework of Theorem~\ref{thmalex} too, but
in order not to overextend the exposition we do not analyse it here in
a very detailed way.

Another classical approach of exponential ergodicity is the spectral
approach, that is, the existence of a \emph{spectral gap}. Recall the
well known equivalence

\begin{proposition}\label{poinc}
Pick $\theta>0$. For any $f \in\mathbb L^2(\pi)$,
\[
\Var_{\pi}( P_t f) \leq \mathrm{e}^{-\theta t}
\Var_{\pi}(f)
\]
if and only if $\pi$ satisfies the following Poincar\'e inequality
\[
\Var_{\pi}(f) \leq\frac{1}{\theta} \int|\nabla f|^2 \,\mathrm{d}
\pi.
\]
\end{proposition}

When $\pi$ is not only invariant but reversible, it is known that both
approaches coincide, that is, the following theorem holds.

\begin{theorem}\label{thmequiv}
The process $Z$ is exponentially ergodic if and only if $\pi$ satisfies
some Poincar\'e inequality. Furthermore if $\pi$ satisfies a Poincar\'e
inequality, we
may choose $\beta=\theta/2$, while if the process is exponentially
ergodic we may choose $\theta=\beta$.
\end{theorem}

The difficult part of this equivalence (i.e., exponential ergodicity
implies Poincar\'e) is shown in Bakry, Cattiaux and Guillin \cite{BCG},
Theorem~2.1, or Kulik \cite{Kul1},
Theorem~3.4.
The converse direction is explained in Cattiaux, Guillin and Zitt \cite{CGZ}.

As it was pointed out in details in the recent papers Kulik \cite
{Kul3} and
Cattiaux, Guillin and Zitt~\cite{CGZ} (which despite the dates of
publication were achieved
simultaneously),
these two formulations of exponential ergodicity are also equivalent to
the existence of exponential moments for the hitting time of a compact set,
in the case of usual diffusion processes. For the reflected processes
we are looking at, some extra work is necessary.

In this paper, we use the second approach related to Theorem~\ref
{thmequiv} to analyse the 2-ball system in the next section, and the
first method presented
in Theorem~\ref{thmalex} to prove the ergodicity of the 3-ball system
developed in
the third section.
For that system, we prove geometric ergodicity in Theorem~\ref
{th_exp_erg_3}, but during the proof
(see, e.g., the statement of Proposition~\ref{prop_sect_3}) we show
that the total quadratic energy $V$ is a Lyapunov function in the sense of
the function $\Phi$ of Theorem~\ref{thmalex}.

The proof is quite intricate.
The key idea is to study and control the hitting time of a \textit{cluster},
that is, a set of relatively small quadratic energy.
It turns out that the most practical way to describe the triangle
configuration built by the three centers is to look at the medians of this
triangle. The
reason is that one has to control a single local time term.

We expect that for any $n=4,5,\ldots$ the system of $n$ stochastic hard
balls will exhibit the same principal behavior: the hitting time of a
(properly defined)
cluster should verify an analogue of Proposition~\ref{prop_sect_3}
which would yield an exponential convergence rate to the invariant measure.
However, the proof for $n=3$ uses specific and comparatively simple
geometry of a 3-ball system: essentially, there exists only one type of
non-clustered
configuration which is \textit{bad} in the sense that some collision could
happen and increase a local time term: two balls are close, while one is
distant. For greater $n$ this method does not work straightforwardly,
and one should take into account the more complicated structure of {\it
sub-clusters}
of close balls: The analysis of local time terms, generated by the
collisions of the balls in this sub-cluster, is much more
delicate.
This is a subject of our further research, and we plan to control the
impact of sub-clusters using induction by $n$. Note that the proof of
Proposition~\ref{prop_sect_3} implicitly contains the induction
step from $n=2$ to $n=3$.

Throughout the proofs, we have tried to trace the constants as
precisely as possible,
in particular to obtain the convergence rate as an explicit function of
$a$. This goal was achieved completely in the case $n=2$ and partially
in the case
$n=3$, where we give explicit estimates on exponential moments of
hitting times of clusters. Clearly, these estimates then would yield a
bound for the
convergence rate, but we do not give it explicitly because of the lack
of explicit formula for $\beta$ in Theorem~\ref{thmalex}.

\section{The case of two balls}\label{sec2}

In this section, we consider the ``baby model'' case $n=2$.
The relative position of the two balls is described by the $\R
^d$-valued process $Y:=Y_1= \frac{X_1-X_2}{2}$ which satisfies
\[
Y(t)=Y(0) + B_1(t) - 2a \int_0^t
Y(s) \,\mathrm{d}s + 2 \int_0^t Y(s) \,\mathrm{d}L_1(s)
,
\]
where $B_1$ is a Brownian motion with covariance $(1/2) \mathrm{Id}$ and
$L_1$ is the local time of $Y$ on the centered sphere with radius
$r/2$, that is, $Y$ is simply
an Ornstein--Uhlenbeck process outside the ball of radius $r/2$ and
normally reflected on the boundary of this ball.
In particular
\[
\pi_a(\mathrm{d}y)= Z_a^{-1} \mathrm{e}^{-4a |y|^2}
\un_{|y|>r/2} \,\mathrm{d}y
\]
is simply a centered Gaussian measure restricted to $D'=\R^d - B(0,r/2)$.

This measure is thus spherically symmetric and radially log-concave, so
that one can use the method in Bobkov \cite{Bobsphere} in order to evaluate
the Poincar\'e
constant.
The following result is a direct consequence of Boissard \textit{et
al.} \cite{BCGM}

\begin{theorem}\label{thm2}
$\pi_a$ satisfies a Poincar\'e inequality with constant $C_P(\pi
_a)=1/\theta_a$ satisfying
\[
\frac{1}2 \biggl(\frac{1}{8a} + \frac{r^2}{4d} \biggr) \leq\max
\biggl(\frac{1}{8a} , \frac{r^2}{4d} \biggr)\leq C_P(
\pi_a) \leq\frac{1}{4a} + \frac{r^2}{4d}.
\]
\end{theorem}

As we explained before, this result captures both the rate of
convergence to a ``well packed'' configuration and the rate of
stabilization of an uniform
rotation.
If we want to avoid the last property, we are led to look at the $\R
^+$-valued process $y(t):=|Y(t)|$ which is the \emph{radial Ornstein--Uhlenbeck}
process reflected at $r/2$ that is, solves the following one
dimensional (at $r/2$) reflected SDE
%
\begin{equation}
\label{eqradial} \mathrm{d}y(t) = \frac{1}{\sqrt{2}} \,\mathrm{d}W(t) - 2a y(t) \,\mathrm{d}t +
\frac
{d-1}{4y(t)}
\,\mathrm{d}t + 2 y(t) \,\mathrm{d}L(t) ,
\end{equation}
with a standard Brownian motion $W$. Its one dimensional reversible
probability measure~is
%
\begin{equation}
\label{eqrev1d} \nu_a(\mathrm{d}\rho) = Z_a^{-1}
\rho^{d-1} \mathrm{e}^{-4a \rho^2} \un_{\rho
>r/2} \,\mathrm{d}\rho,
\end{equation}
for which we have the following result which furnishes a bound of the
rate of ``packing'' of two balls.

\begin{theorem}\label{thm2bis}
Let $(P_t)_{t\ge0}$ be the transition distribution of the half distance
$(y(t))_{t\ge0}$ between the centers of the two balls moving according
to the dynamics
\emph{(A)}.
\[
\forall y>\frac{r}{2} \qquad\bigl\Vert P_t(y,\cdot) - \nu_a
\bigr\Vert_{\mathrm{TV}} \le C(y) \mathrm{e}^{-4a t}.
\]
\end{theorem}

\begin{pf}
Theorem~\ref{thm2bis} holds as soon as $\nu_a$ satisfies a Poincar\'e
inequality with constant $C_P(\nu_a)$ satisfying
$C_P(\nu_a) \le\frac{1}{8a}$.
This latter result is a simple application of Bakry--Emery criterion on
the interval $\rho\geq r/2$.
Recall that Bakry--Emery criterion tells us that provided $V''(\rho
)\geq
A>0$ for all $\rho\in\mathbb R$, then the (supposed to be finite)
measure $\mathrm{e}^{-V(\rho)}
\,\mathrm{d}\rho$ satisfies a Poincar\'e inequality with constant $1/A$.
The measure $\mathrm{e}^{-V(\rho)} \un_{\rho>r/2} \,\mathrm{d}\rho$ can be approximated,
as $N \to+\infty$, by $\mathrm{e}^{-(V(\rho)+N((r/2)-\rho)_+^4)} \,\mathrm{d}\rho$
which still
satisfies the same lower bound for the second derivative, uniformly in
$N$, showing that Bakry--Emery criterion extends to the case of an interval.
Finally, it is immediately seen that $\nu_a$ satisfies Bakry--Emery
criterion with $A=8a$.
\end{pf}

\begin{figure}[b]

\includegraphics{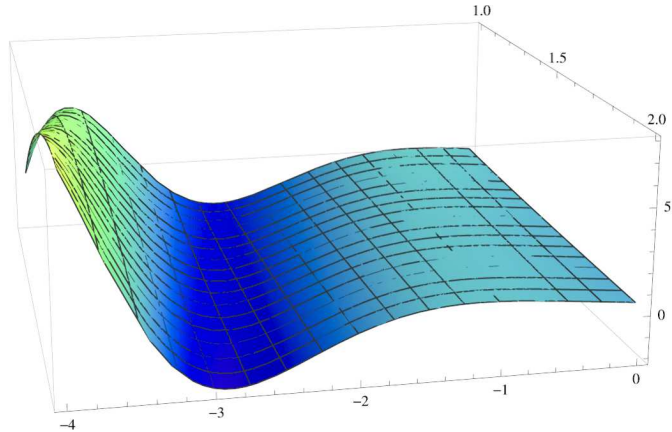}

\caption{3-dimensional plot of the function $(a,\mathrm{x}) \mapsto U(
\mathrm{x} +1, 2, a)$.}\label{fig1}
\end{figure}

\begin{pf*}{Sketch of the proof of Theorem~\ref{thm2}}
Actually $\nu_a$ is the radial part of $\pi_a$. In polar coordinates
$(\rho, s) \in\mathbb R \times S^{d-1}$, $\pi_a$ factorizes as $\nu_a
\otimes {d}s$ where
${d}s$ is the normalized\vadjust{\goodbreak} uniform measure on the sphere $S^{d-1}$, which
satisfies a Poincar\'e inequality with constant $1/d$.
Bobkov's method exposed in Bobkov \cite{Bobsphere} and detailed in
Boissard \textit{et al.} \cite{BCGM},
Proposition~2.1, allows us to deduce the upper bound of Theorem~\ref
{thm2}. For the
lower bound, it is enough to consider linear
functions.
\end{pf*}

\begin{figure}[t]

\includegraphics{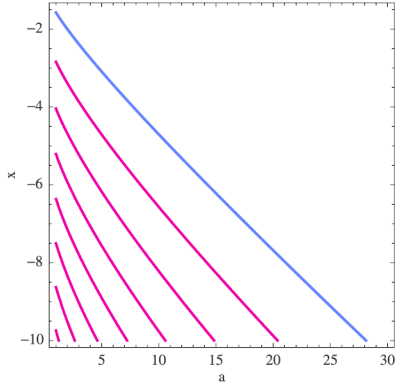}

\caption{The zeros of the function $(a,\mathrm{x}) \mapsto U( \mathrm{x}
+1, 2, a)$.}
\label{fig2}
\end{figure}

\begin{figure}[b]

\includegraphics{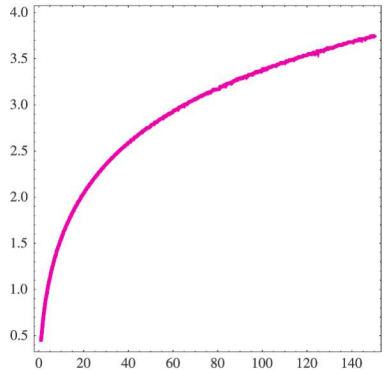}

\caption{The logarithm of $-x_a$ as function of $a$.}
\label{fig3}
\end{figure}

\begin{remark}\label{remlinetsky}
The spectral gap ($\theta/2$ in Proposition~\ref{poinc}) of linear
diffusion processes can be studied by solving some O.D.E.
For instance in our case, if we consider the process
$z(t):=y(t)^2-\frac
{r^2}{4}$, it is an \emph{affine diffusion} reflected at $0$,
that is, it solves the reflected~SDE\looseness=-1
\[
\mathrm{d}z(t)=\sqrt2 \sqrt{z(t)+\bigl(r^2/4\bigr)} \,\mathrm{d}B_t +4a
\biggl(\frac{d}{8a} - \frac{r^2}{4} - z(t) \biggr) \,\mathrm{d}t +
\mathrm{d}L_z(t) ,
\]\looseness=0
where $L_z(\cdot)$ is proportional to the local time of the process
$z$ at $0$.
For such linear processes, it is shown in Linetsky \cite{Linets},
Section~6.2,
that the spectral gap is given by
$\theta/2 = - 4a x $ where $x_a$ is the first negative zero of the
Tricomi confluent hypergeometric function
\[
\mathrm{x} \mapsto U\biggl( \mathrm{x} +1, 1 + \frac{d}{2}, a
r^2\biggr) ,
\]
which is not easy to calculate.
Nevertheless the following Figures \ref{fig1} to \ref{fig3} -- done
by simulations using
Mathematica9 and the built-in functions \emph{FindRoot} and \emph
{HypergeometricU}
for $d=2$ and $r=1$ -- lead to conjecture that $a \mapsto x_a$ is bounded,
and therefore the spectral gap of the process $z$ is sublinear in $a$.
In the Figure~\ref{fig2} the curve, being the
upper most one, corresponds to the function $a \mapsto x_a$. Scrolling
the picture from up to down, one meets the
curve corresponding to the second negative zero and so on.\vadjust{\goodbreak}
\end{remark}

\begin{remark}\label{remlyap}
It is well known that $\Phi(z)=|z|^2$ plays the role of a Lyapunov
function for the Ornstein--Uhlenbeck process.
But for the radial Ornstein--Uhlenbeck process $y(\cdot)$ the
situation is
still better.
Indeed its infinitesimal generator denoted by ${\mathcal L}$ is given,
for $\rho>r/2$, by
\[
{\mathcal L} g(\rho) = \frac{1}4 g''(\rho)
- \biggl(2a \rho- \frac
{d-1}{4\rho} \biggr) g'(\rho) ,
\]
so that, if $g(\rho)=\phi(\rho)=\rho^2$, it holds
\[
{\mathcal L} \phi(\rho)=\frac{d}2 - 4a \rho^2\qquad
\mbox{provided $\rho>r/2$} ,
\]
so that ${\mathcal L}\phi\leq-4a\varepsilon\phi$ as soon as $\frac
{1}{1-\varepsilon} \leq\frac{2 ar^2 }{d}$ for some $\varepsilon>0$.

It follows that the process $t \mapsto \mathrm{e}^{4a \varepsilon t} y^2(t)$ is
a supermartingale up to the first time $\tau_r$ the process $y$ hits
the value $r/2$.

This yields
the following proposition.
\end{remark}

%
\begin{proposition}\label{prophit}
Assume that $a>\frac{d}{2r^2}$ and define $\tau_r=\inf\{t ;
y(t)=|Y(t)|=r/2\}$ the hitting time of a packing configuration. Then
\[
\mathbb\Ps(\tau_r>t) \leq\frac{4 \Es(|Y(0)|^2)}{r^2} \mathrm{e}^{-
4(a-{d}/{(2r^2)}) t}.
\]
\end{proposition}

This means that
the system reaches a packing configuration $y=|Y|=r/2$ before time $t$
with a probability at least equal to
$1-\frac{4 \Es(|Y(0)|^2)}{r^2} \mathrm{e}^{- 4(a-{d}/{(2r^2)})t}$.
This statement should be particularly interesting to generalize to a
higher number of balls.

In the next section, we shall look at the case of three hard balls,
where real difficulties begin to occur.

\section{The case of three balls}\label{sec3}

We now address the case of three Brownian hard balls with attractive
interaction.
That is, we consider the dynamics (A) and (B) with $n=3$ and $a>0$. For
simplicity, from now on we assume that $r=1$.

Our aim is to prove Theorem~\ref{thmmain} for $n=3$. Recall that the
relative positions
\[
Y_i=X_i-(X_1+X_2+X_3)/3
\]
of the three hard balls follow the dynamics
\[
(\mathrm{B}) \cases{ %
\mbox{for } i \in\{1,2,3\} , t \in
\R^+ ,
\vspace*{2pt}\cr
\displaystyle Y_i(t) = Y_i(0) + W_i(t)-
\overline{W}(t) - a \sum_{j=1}^3 \int
_0^t \bigl(Y_i(s)-Y_j(s)
\bigr) \,\mathrm{d}s \vspace*{2pt}\cr
\displaystyle \hspace*{32pt}{}+ \sum_{j=1}^3 \int
_0^t \bigl(Y_i(s)-Y_j(s)
\bigr) \,\mathrm{d}L_{ij}(s),
\vspace*{2pt}\cr
L_{ij}(0) = 0 , L_{ij} \equiv L_{ji}
\mbox{ and }\displaystyle L_{ij}(t) = \int_0^t
\un_{|Y_i(s)-Y_j(s)|=1} \,\mathrm{d}L_{ij}(s), L_{ii} \equiv0,}
\]
where $W_1$, $W_2$, $W_3$ are independent standard $d$-dimensional
Brownian motions and $\overline{W}:=(W_1+W_2+W_3)/3$.
We noted in Section~\ref{Intro} that the process $Y$ is a $\bar
{D}'$-valued Feller process satisfying the local Dobrushin condition, where
\[
D' = \bigl\{\by\in\bigl(\mathbb R^d
\bigr)^3 ; |y_1 - y_2|>1 , |y_2 -
y_3|>1 , |y_3 - y_1|>1 \mbox{ and }
y_1+y_2+y_3=0 \bigr\}.
\]
The invariant probability measure for the system (B) is given by
\[
\mathrm{d}\pi_a(\by)= Z_a^{-1} \mathrm{e}^{- ({a}/{2}) V(\by)}
\un_{
D'}(\by) \,\mathrm{d}\by,
\]
where $Z_a$ is the normalization constant. The function $V(\by
)=|y_1-y_2|^2+|y_2-y_3|^2+|y_3-y_1|^2$ is, as before, the quadratic
energy of the system.
The configurations with minimal quadratic energy are triangular
packings and build the set
\[
\EE_{\mathrm{min}} = \bigl\{ \by\in D' ; V(\by)=3 \bigr\}.
\]
We will prove in the sequel the following theorem.

\begin{theorem}\label{th_exp_erg_3} 
Let $Y$ satisfying \emph{(B)} be the process of the relative positions
of three hard balls and let $(P_t)_t$ be its transition distribution.
There exists
$\beta>0$ such that
\[
\forall\by\in D' \qquad \bigl\Vert P_t(\by,\cdot) -
\pi_a \bigr\Vert_{\mathrm{TV}} \le C(\by) \mathrm{e}^{-\beta t}.
\]
In the large attraction regime, we obtain asymptotically in time the
concentration of the system around packing triangular configurations in
the sense that,
for all $\by\in D'$,
%
\begin{eqnarray}
\label{asymptotic_a} &&\forall\eps,\eta>0, \exists a_0, t_0\quad
\mbox{s.t.}\quad a>a_0 \mbox{ and}
\nonumber
\\[-8pt]
\\[-8pt]
\nonumber
&&\quad t>t_0\quad \Rightarrow\quad \Ps
\bigl(\operatorname{dist}\bigl(Y(t),\EE_{\mathrm{min}}\bigr)\le\eta|
Y(0)=\by
\bigr) \ge1-\eps.
\end{eqnarray}
\end{theorem}


\begin{pf*}{Proof of (\ref{asymptotic_a})}
Let $\EE_{\mathrm{min}}^\eta= \{ \by\in D' ; V(\by) \le3+\eta\}$ be the
set of configurations with $\eta$-minimal energy.
Clearly, $\pi_a(\EE_{\mathrm{min}}^\eta)$ is large for $a$ large enough:
\[
\forall\eps>0\ \exists a_0 \mbox{ s.t. } \forall a>a_0\qquad
\pi_a\bigl(\EE_{\mathrm{min}}^\eta\bigr) \ge1-\eps
\]
because
\[
\pi_a\bigl(\bigl(\EE_{\mathrm{min}}^\eta
\bigr)^c\bigr) =\frac{ \int_{D'} \mathrm{e}^{- ({a}/{2}) V(\by)} \un
_{V(\by)>3+\eta
} \,\mathrm{d}\by}{ \int_{D'} \mathrm{e}^{- ({a}/{2}) V(\by)} \,\mathrm{d}\by}
\le\frac{ \int
_{D'} \mathrm{e}^{- ({a}/{2}) (V(\by)-3-\eta)} \un
_{V(\by)>3+\eta} \,\mathrm{d}\by}{ \int_{D'} \mathrm{e}^{- ({a}/{2}) (V(\by)-3-\eta
)} \un_{V(\by) \le
3+\eta} \,\mathrm{d}\by}
\]
hence
\[
\pi_a\bigl(\bigl(\EE_{\mathrm{min}}^\eta
\bigr)^c\bigr) \le\frac{ 1 }{ \int_{D'} \un_{V(\by) \le3+\eta}
\,\mathrm{d}\by} \int_{D'}
\mathrm{e}^{-({a}/{2}) (V(\by)-3-\eta)} \un_{V(\by)>3+\eta} \,\mathrm{d}\by,
\]
which vanishes for $a$ tending to infinity. On the other side, the
convergence in total variation implies
\[
\lim_{t \to+\infty} \Ps\bigl( Y(t)\in\EE_{\mathrm{min}}^\eta
|Y(0)=\by\bigr)=\pi_a\bigl(\EE_{\mathrm{min}}^\eta
\bigr).
\]
Using the continuity of $V$, this yields (\ref{asymptotic_a}).

The same technique obviously works for any number $n$ of balls.
\end{pf*}

The technique we will use in order to prove the main part of the above
theorem, that is, the exponential ergodicity, is very intricate. It relies
on hitting time
estimates and is the subject of the rest of the paper.

\subsection{Quadratic energy and hitting time of clusters}\label
{subsec31}

In this Section~\ref{subsec31}, we present the energy as a Lyapounov
function, define the compact set of cluster patterns and state that $Y$
satisfies the
assumptions of Theorem~\ref{thmalex}.

For a configuration $\bx=(x_1,x_2,x_3) \in\R^{3d}$, or equivalently
for a pattern $\by=(y_1,y_2,y_3) \in\R^{3d}$
of relative positions (with $y_i=x_i-(x_1+x_2+x_3)/3$), recall that the
quadratic (total) energy satisfies
\[
V(\by)= 3 \bigl(|y_1|^2 + |y_2|^2
+ |y_3|^2\bigr) = |x_1-x_2|^2
+ |x_2-x_3|^2 + |x_3-x_1|^2.
\]

\begin{definition}\label{defRcluster}
Fix $R>0$. We say that a relative position $\by=(y_1,y_2,y_3)\in\R
^{3d}$ forms an \emph{$R$-cluster},
if there exists a permutation $\sigma$ on $\{1,2,3\}$ such that
\[
|y_{\sigma(1)}-y_{\sigma(2)}|^2 \leq1+R\quad \mbox{and}\quad
|y_{\sigma(2)}-y_{\sigma(3)}|^2 \leq1+R,
\]
or equivalently,
\[
\forall i \in\{1,2,3\}, \exists j \neq i,\qquad |y_i-y_j|^2
\leq1+R.
\]
\end{definition}

Note that, since $y_1+y_2+y_3=0$,
\[
K_R := \bigl\{\by\in\bigl(\R^d\bigr)^3 ; \by
\mbox{ forms an $R$-cluster} \bigr\}
\]
is a compact set of $\R^{3d}$.

In order to check the assumptions of Theorem~\ref{thmalex}, we would
like to control the time needed by the process $Y$
in such a way that its value forms an $R$-cluster.

\begin{proposition}\label{prop_sect_3}
The relative positions process $Y$ has the following properties:
\begin{longlist}[1.]
\item[1.] Its Lyapounov quadratic energy $V$ fulfills
%
\begin{equation}
\label{verif1} \forall\by\in D'\qquad V(\by) \ge3\quad \mbox{and} \quad\lim
_{|\by
|\to+\infty}V(\by)=+\infty.
\end{equation}
\item[2.] The $K_R$-hitting time
\[
\tau:= \inf\bigl\{t\geq0 , Y(t)=\bigl(Y_1(t),Y_2(t),Y_3(t)
\bigr) \mbox{ forms an } R \mbox{-cluster}\bigr\}
\]
satisfies inequalities
%
\begin{equation}
\label{1} \forall\by\in D'\qquad \Es_\by\bigl(
\mathrm{e}^{\lambda\tau} V\bigl(Y(\tau)\bigr) \bigr) \leq V(\by)
\end{equation}
and
%
\begin{equation}
\label{ineq_energy_clusters} \forall\by\in D'\qquad \Es_\by\bigl( V
\bigl(Y(t)\bigr) \un_{\tau>t} \bigr) \le2 \mathrm{e}^{-\lambda t} V(\by)
\end{equation}
for
%
\begin{equation}
\label{hyp_lambda_R} \lambda=\min\bigl(a,a^2\bigr) \quad\mbox{and any}\quad R \ge
\frac
{48a+16d+60}{a} \mathrm{e}^{10\,505/a}.
\end{equation}

\item[3.]
The quadratic energy of the system is uniformly bounded in time for
any initial $R$-cluster position, that is, for $R$ as above
%
\begin{equation}
\label{energy_clusters_5} \sup_{\by\in K_R} \sup_{t>0}
\Es_\by\bigl( V\bigl(Y(t)\bigr) \bigr) <+\infty.
\end{equation}
\end{longlist}
\end{proposition}


\begin{pf*}{Proof of Proposition~\ref{prop_sect_3}}
The properties (\ref{verif1}) are an obvious consequence of the
definition of~$V$.

The proofs of (\ref{1}) and (\ref{ineq_energy_clusters}) rely on a
study of the length of the largest median in the triangle of particles,
as a function of the
time. These proofs are postponed to Section~\ref{subsec32}.

In order to prove (\ref{energy_clusters_5}), we first establish the
following consequence of (\ref{ineq_energy_clusters})
%
\begin{equation}
\label{energy_clusters_28} \exists T>0\ \exists D_1>0 \quad\mbox{s.t.}\quad \sup
_{\by\in K_R} \sup_{t\in[0;T]} \Es_\by\bigl(
V\bigl(Y(t)\bigr) \bigr) \le D_1.
\end{equation}
Take $R$ as in Proposition~\ref{prop_sect_3} part $(2)$ and take some
larger $\bar{R}>R$.
Construct a sequence of stopping times $\xi_k, k\geq0$ in the
following way. Assume that $Y(0)=\by\in K_R$ and define $\xi_0=0$,
\[
\xi_{2j-1}=\inf\bigl\{t>\xi_{2j-2}\dvt Y(t)\notin
K_{\bar{R}}\bigr\},\qquad \xi_{2j}=\inf\bigl\{t>\xi_{2j-1}
\dvt Y(t)\in K_{R}\bigr\},\qquad j\geq1.
\]
Then
\[
\Es_\by\bigl( V\bigl(Y(t)\bigr) \bigr) =\sum
_{k=1}^\infty\Es_\by\bigl( V\bigl(Y(t)
\bigr)\un_{t\in[\xi_{k-1}, \xi_k)} \bigr) =:\sum_{k\mathrm{\ is\ even}}+\sum
_{k\mathrm{\ is\ odd}}.
\]
Let
\[
\| V \|_{K_{\bar{R}}}:=\max_{ \by\in K_{\bar{R}} } V(\by)=6(1+\bar{R}).
\]
When $k$ is odd and $t\in[\xi_{k-1}, \xi_k)$, we have $Y(t)\in
K_{\bar
{R}}$, which means that
\[
\sum_{k\mathrm{\ is\ odd}} \le\| V \|_{K_{\bar{R}}}.
\]
When $k$ is even, we have
\[
\Es_\by\bigl( V\bigl(Y(t)\bigr)\un_{t\in[\xi_{k-1}, \xi_k)}
\bigr) \leq
\Es_\by\bigl( \un_{t\geq\xi_{k-1}} \bigl(\Es\bigl[V\bigl
(Y(t)\bigr)\un
_{\xi
_{k}>t}|\Ff_{\xi_{k-1}}\bigr] \bigr) \bigr).
\]
Note that, by the continuity of trajectories, $Y(\xi_{k-1})\in K_{\bar{R}}$.
Then, applying the strong Markov property at the time moment $\xi
_{k-1}$ and
(\ref{ineq_energy_clusters}), we get
\[
\Es_\by\bigl( V\bigl(Y(t)\bigr)\un_{t\in[\xi_{k-1}, \xi_k)}
\bigr) \le2 \| V
\| _{K_{\bar{R}}} \Ps_\by(t\geq\xi_{k-1}).
\]
Hence,
\[
\sum_{k\mathrm{\ is\ even}}\leq2 \| V \|_{K_{\bar{R}}} \sum
_{k\mathrm{\ is\ even}}\Ps_\by(t\geq\xi_{k-1}),
\]
and to prove (\ref{energy_clusters_28}) it is enough to prove that for
some $T$
\[
\sup_{\by\in K_R} \sup_{t\le T}\sum
_{k\mathrm{\ is\ even}}\Ps_\by(t\geq\xi_{k-1})<\infty.
\]
By the Chebyshev--Markov inequality, for any fixed $c>0$ and $T>0$
\[
\Ps_\by(t\geq\xi_{k-1}) \le \mathrm{e}^{ct}
\Es_\by\bigl( \mathrm{e}^{-c\xi_{k-1}} \bigr) \le \mathrm{e}^{cT}
\Es_\by\bigl( \mathrm{e}^{-c\xi_{k-1}} \bigr)\qquad \forall t\leq T.
\]
Clearly, the exponential moment $\Es_\by(\mathrm{e}^{-c\xi_{k-1}})$ can be
expressed iteratively via the conditional exponential moments of the differences
$\xi_j-\xi_{j-1}$ w.r.t. $\Ff_{\xi_{j-1}}, j=1, \ldots, k-1$. When $j$
is odd, this conditional exponential moment can be estimated as follows:
\[
\Es_\by\bigl[\mathrm{e}^{-c(\xi_j-\xi_{j-1})} |\Ff_{\xi_{j-1}} \bigr]
\le\sup
_{\by\in K_R}\Es_\by\bigl(\mathrm{e}^{-c\varsigma}\bigr),\qquad
\varsigma=\inf\bigl\{t\dvt X(t)\notin K_{\bar{R}}\bigr\}.
\]
Note that
\[
q:=\sup_{\by\in K_R}\Es_\by\bigl(\mathrm{e}^{-c\varsigma}
\bigr)<1
\]
because otherwise, by the Feller property of the process $Y$, there
would exist $\by\in K_R$ such that $\varsigma=0$ $\Ps_\by$-a.s., which
would contradict the
continuity of the trajectories of $Y$. Then
\begin{eqnarray*}
\Es_\by\bigl(\mathrm{e}^{-c\xi_{k-1}}\bigr) &\le& q^{k/2},
\\
\sup_{\by\in K_R} \sup_{t\le T}\sum
_{k\mathrm{\ is\ even}} \Ps_\by(t\ge\xi_{k-1})&\le&
\mathrm{e}^{cT}\sum_{k\mathrm{\ is\ even}}q^{k/2}<\infty,
\end{eqnarray*}
which completes the proof of (\ref{energy_clusters_28}).

Let us now deduce the uniform bound (\ref{energy_clusters_5}) from the
finite time bound (\ref{energy_clusters_28}).
This proof is simple and similar to that of
Lemma A.4 in Kulik \cite{Kul1}. Indeed, let $\tau'$ be the first time moment
for $X(t)$ to form an $R$-cluster after $T$
\[
\tau':=\inf\bigl\{ t \ge T \mbox{ s.t. } Y(t) \mbox{ forms an
$R$-cluster} \bigr\}.
\]
Then for $t>T$
\[
\Es_\by\bigl( V\bigl(Y(t)\bigr) \bigr) = \Es_\by\bigl(
V\bigl(Y(t)\bigr)\un_{\tau'>t} \bigr) + \Es_\by\bigl( V
\bigl(Y(t)\bigr)\un_{\tau'\le t} \bigr).
\]
We have by the Markov property of $Y$, for $\lambda$ small enough and
$R$ large enough for (\ref{ineq_energy_clusters}) to hold
\[
\Es_\by\bigl( V\bigl(Y(t)\bigr)\un_{\tau'>t} \bigr) = \int
_{\R^{3d}} \bigl( \Es_\bx\bigl( V\bigl(Y(t-T)\bigr)
\un_{\tau>t-T} \bigr) \bigr)\Ps_T(\by, \mathrm{d}\bx)
\le2\mathrm{e}^{-\lambda(t-T)} \Es_\by\bigl( V\bigl(Y(T)\bigr) \bigr).
\]
On the other hand, by the strong Markov property of $Y$, we have
\[
\Es_\by\bigl( V\bigl(Y(t)\bigr)\un_{\tau'\le t} \bigr) =
\Es_\by\bigl( \bigl( \Es_{Y(\tau')}V\bigl(Y\bigl(t-
\tau'\bigr)\bigr) \bigr) \un_{\tau'\le t} \bigr) \le\sup
_{\by\in K_R} \sup_{s\le t-T} \Es_\by\bigl(
V\bigl(Y(s)\bigr) \bigr).
\]
Let $D_k=\sup_{\by\in K_R}\sup_{t\le kT}\Es_\by( V(Y(t))
)$.
Then the above estimates and (\ref{energy_clusters_28}) yield for every
$\by\in K_R$ and
$(k-1)T \le t\le kT$
\[
\Es_\by\bigl( V\bigl(Y(t)\bigr) \bigr) \le2\mathrm{e}^{-\lambda
(k-2)T}\sup
_{\by
\in
K_R}\Es_\by\bigl(V\bigl(Y(T)\bigr)
\bigr)+D_{k-1} \le2\mathrm{e}^{-\lambda(k-2)T}D_1+D_{k-1}.
\]
Then
\[
D_k=\max\Bigl(D_{k-1}, \sup_{\by\in K_R}
\sup_{(k-1)T \le t\le kT} \Es_\by\bigl( V\bigl(Y(t)\bigr) \bigr)
\Bigr) \le2\mathrm{e}^{-\lambda(k-2)T}D_1+D_{k-1},
\]
and consequently
\[
\sup_{\by\in K_R} \sup_{t\ge0} \Es_\by
\bigl( V\bigl(Y(t)\bigr) \bigr) \le D_1 +D_1\sum
_{k=2}2\mathrm{e}^{-\lambda(k-2)T}= D_1+2D_1
\bigl[1-\mathrm{e}^{-\lambda
T}\bigr]^{-1}<\infty.
\]
\upqed\end{pf*}
%

\subsection{Cluster hitting time estimates}\label{subsec32}

This section is devoted to the proof of assertion $(2)$ of Proposition~\ref{prop_sect_3}, that is (\ref{1}) and (\ref{ineq_energy_clusters}).
It will complete
the proof of Theorem~\ref{prop_sect_3}. From now on, $R$ and $\lambda$
are fixed parameters. If the starting configuration $Y(0)=\by\in\R
^{3d}$ already forms
an $R$-cluster, $\tau=0$ and (\ref{1}), (\ref{ineq_energy_clusters})
are trivial. In the sequel, we assume that $\by\notin K_R$.

Let us reduce the problem to the study of the time the farthest ball
need to come closer to the others. Since $Y(0)=\by\notin K_R$, some
ball $\mathbf{i}$ has a
center $y_i$ which is farther than $\sqrt{1+R}$ from the other two
centers. Suppose, for instance, $\mathbf{i}=\mathbf{1}$. We construct a
sequence of
stopping times
corresponding to hitting times of levels for the distance between the
balls \textbf{2} and \textbf{3}. These levels depend on two
parameters $\medp
>\medm>0$ which
will be chosen later.
Put $\sigma_0=0$, and for $k\geq1$,
\begin{eqnarray*}
\sigma_{2k-1}&:=&\inf\bigl\{ t>
\sigma_{2k-2}\dvt|Y_2-Y_3|^2\le1+2
\medm\bigr\},
\\
\sigma_{2k}&: =&\inf\bigl\{ t>\sigma_{2k-1}\dvt
|Y_2-Y_3|^2 \ge1+2\medp\bigr\},
\end{eqnarray*}
and define a \emph{time-depending border level} as
\[
R(t)=\sum_{k=1}^\infty\bigl[R
\un_{t\in[\sigma_{2k-2}, \sigma
_{2k-1})}+R'\un_{t\in[\sigma_{2k-1}, \sigma_{2k})} \bigr],\qquad t\in
[0,\infty)
\]
for some $0<R'<R$. Let us now define the first time the ball number
\textbf{1} comes closer than $R(\cdot)$ to one of the others
\[
\tau_1 := \inf\bigl\{t\geq0 ; \min\bigl(\bigl|Y_1(t)-
Y_2(t)\bigr|^2 , \bigl|Y_1(t)-Y_3(t)\bigr|^2
\bigr) \le1+ R(t)\bigr\}.
\]
Consider the configuration $Y(\tau_1)$. If it forms an $R$-cluster,
then put $\tau_2=\tau_1$.
Otherwise, there exists some ball whose center is farther than $\sqrt
{1+R}$ from the other two centers.
Since $R(t)\le R$, this ball should be either \textbf{2} or \textbf{3}.
Define then the new sequence of stopping times corresponding to level
hitting times of the distance between the other two balls (with the
same values $\medp$,
$\medm$) and the corresponding time-depending border level and
respective $\tau_2$, and so on.
By monotonicity, there exists an a.s. limit, $\tau_\infty=\lim_n\tau_n$.

To obtain the desired estimates on the $K_R$-hitting time for $Y$, we
only have to prove the following proposition.

\begin{proposition}\label{prop_energy_clusters}
If $R$ and $\lambda$ are as in (\ref{hyp_lambda_R}), for any starting
configuration $\by\in D'$
\[
\Es_\by\bigl( \mathrm{e}^{\lambda\tau_1}V\bigl(Y(\tau_1)\bigr)
\bigr) \le V(\by)
\]
and for any $\by\in D'$ and any finite time horizon $T\in\R^+$
\[
\Es_\by\bigl( \mathrm{e}^{\lambda\tau_1}V\bigl(Y(\tau_1)\bigr) +
\mathrm{e}^{\lambda(\tau_1
\wedge T)}V\bigl(Y(\tau_1\wedge T)\bigr) \bigr) \le2 V(\by).
\]
\end{proposition}

Remark that Proposition~\ref{prop_energy_clusters} implies
Proposition~\ref{prop_sect_3} part $(2)$.

Indeed, by construction and by the strong Markov property of $Y(t),
t\geq0$, it follows from Proposition~\ref{prop_energy_clusters} and
its analogous for
$\tau_2-\tau_1,\tau_3-\tau_2, \ldots,$ that
%
\begin{equation}
\label{3} \Es_\by\bigl( \mathrm{e}^{\lambda\tau_n}V\bigl(Y(
\tau_n)\bigr) \bigr)\leq V(\by),\qquad n\geq1.
\end{equation}
Because $V$ is bounded from below, this implies that $\tau_\infty
<\infty
$ a.s. On the other hand, by the construction and by the continuity of
the trajectories
of $Y(t), t\geq0$ it is easy to see that $Y(\tau_\infty)$ forms an
$R$-cluster as soon as $\tau_\infty<\infty$, and consequently $\tau
\leq
\tau_\infty$ a.s.
Hence, (\ref{1}) follows from~(\ref{3}) by the Fatou lemma.

By Fatou lemma again, the second inequality in Proposition~\ref
{prop_energy_clusters} implies
\[
\Es_\by\bigl(\mathrm{e}^{\lambda(\tau_\infty\wedge T)} V\bigl(Y(\tau
_\infty\wedge
T)\bigr) \bigr) \le2 V(\by).
\]
Since $\tau_\infty\wedge T \ge\tau\wedge T$ this gives
\[
\Es_\by\bigl( \mathrm{e}^{\lambda(\tau\wedge T)} V\bigl(Y(\tau_\infty
\wedge
T)\bigr) \bigr) \le2 V(\by),
\]
which implies (\ref{ineq_energy_clusters}) because $\mathrm{e}^{\lambda T} \un
_{\tau>T} \le \mathrm{e}^{\lambda(\tau\wedge T)}$ and $\un_{\tau>T}
V(Y(\tau
_\infty\wedge T)) =
\un_{\tau>T} V(Y(T))$.

Our aim now is to prove Proposition~\ref{prop_energy_clusters}.

\subsection{Dynamics of the medians of the triangle}\label
{subsec_median}

Let us introduce the following vectors describing the triangle $Y_1,Y_2,Y_3$:
\[
U_1:=\sqrt{\frac{2}{3}} \biggl(\frac{Y_2+Y_3}{2}-Y_1
\biggr),\qquad  U_{23}:=\frac{1}{\sqrt2}(Y_2-Y_3).
\]
$U_1$ is the (scaled) median starting from $Y_1$ and $U_{23}$ is its
(scaled) opposite side.

In order to prove Proposition~\ref{prop_energy_clusters}, we only have
to consider the behaviour of $Y$ up to time $\tau_1$.
But, before the time moment $\tau_1$, the ball \textbf{1} does not
hit any
other ball. Therefore, on the (random) time interval $[0, \tau_1]$
processes $U_1$,
$U_{23}$ satisfy the following simple SDEs:
%
\begin{equation}
\label{4} \cases{ %
\mathrm{d}U_1(t)=\mathrm{d}B_1(t)-3aU_1(t)
\,\mathrm{d}t,
\vspace*{2pt}\cr
\mathrm{d}U_{23}(t)=\mathrm{d}B_{23}(t)-3aU_{23}(t)
\,\mathrm{d}t+2U_{23}(t) \,\mathrm{d}L_{23}(t).}
\end{equation}
Note that the martingale terms $B_1$, $B_{23}$ are independent $\R
^d$-valued Brownian motions and that the dynamics of the median $U_1$
does not include a local
time term up to time $\tau_1$.

Also note that the quadratic energy has a simple expression as a
function of the median and its opposite side, and that these two
lengths control the size of
the triangle $Y$.

\begin{lemma}
\[
V(Y)=3\bigl(|U_1|^2+|U_{23}|^2
\bigr)
\]
and for $j=2$ or $j=3$
%
\begin{equation}
\tfrac{1}{3}|Y_2-Y_1|^2+
\tfrac{1}{3}|Y_3-Y_1|^2-
\tfrac{1}{3}|U_{23}|^2 = |U_1|^2
\le\tfrac{4}{3}|Y_j-Y_1|^2+
\tfrac{2}{3}|U_{23}|^2. \label
{eq:boundsmedians}
\end{equation}
\end{lemma}
%

\begin{pf}
The equalities are simple norm computations and $U_1=\frac{1}{\sqrt
{6}}((Y_3-Y_2)-2(Y_1-Y_2))$, that is, $|U_1|^2 =\frac
{1}{6}|2(Y_2-Y_1)-\sqrt
{2}U_{23}|^2$ gives
the upper bound for $j=2$.
\end{pf}

\subsection{The time weighted energy decreases}\label
{subsec_weighted_energy}

Define the time weighted energy of the system by
\[
H(t):=\mathrm{e}^{\lambda t}V\bigl(Y(t)\bigr) =3\mathrm{e}^{\lambda t}\bigl|U_1(t)\bigr|^2
+ 3\mathrm{e}^{\lambda
t}\bigl|U_{23}(t)\bigr|^2
\]
and compute its mean value $\Es_\by( H(\zeta_k) )$ at the
random time $\zeta_k:=T \wedge\tau_1\wedge\sigma_k$ for a fixed time
horizon $T$ which
may be finite ($T\in\R^+$ and $\zeta_k:=T \wedge\tau_1\wedge
\sigma
_k$) or infinite ($T=+\infty$ and $\zeta_k:=\tau_1\wedge\sigma_k$).
%
\begin{eqnarray}
\label{5} &&\Es_\by\bigl( H(T \wedge\tau_1\wedge
\sigma_k) + H(\tau_1\wedge\sigma_k)
\bigr)
\nonumber
\\
&&\quad=\Es_\by\bigl( \bigl(H(T \wedge\tau_1\wedge
\sigma_k) + H(\tau_1\wedge\sigma_k)\bigr)
\un_{T \wedge\tau_1 > \sigma_{k-1}} \bigr)
\nonumber
\\[-8pt]
\\[-8pt]
\nonumber
&&\qquad{} +\Es_\by\bigl( \bigl(H(T \wedge\tau_1\wedge
\sigma_{k-1}) + H(\tau_1\wedge\sigma_k)\bigr)
\un_{\tau_1 > \sigma_{k-1} \ge T} \bigr)
\\
&&\qquad{} +\Es_\by\bigl( \bigl(H(T \wedge\tau_1\wedge
\sigma_{k-1}) + H(\tau_1\wedge\sigma_{k-1})
\bigr) \un_{\sigma_{k-1} \ge\tau_1}\bigr) .
\nonumber
\end{eqnarray}
But
\[
\Es_\by\bigl( H(T \wedge\tau_1\wedge
\sigma_k)\un_{T \wedge\tau
_1>\sigma_{k-1}} \bigr) =\Es_\by\bigl(
\un_{T \wedge\tau_1> \sigma_{k-1}} \mathrm{e}^{\lambda
\sigma
_{k-1}} \Es_\by\bigl[H^{\sigma_{k-1}}(T
\wedge\tau_1\wedge\sigma_k) |\Ff_{\sigma_{k-1}} \bigr]
\bigr),
\]
where
\[
H^\sigma(t):=\mathrm{e}^{\lambda( t-\sigma\wedge t)}V\bigl(Y(t)\bigr).
\]

\begin{proposition}\label{lem_recur}
Under a proper choice of $\lambda$, $R'$ and $R$, for every $k \ge1$
%
\begin{equation}
\label{6} \Es_\by\bigl[H^{\sigma_{k-1}}(\tau_1
\wedge\sigma_k) |\Ff_{\sigma
_{k-1}} \bigr] \le H^{\sigma_{k-1}}(
\sigma_{k-1}) \qquad\mbox{on the set } \{\tau_1>
\sigma_{k-1}\}
\end{equation}
and for each finite time horizon $T\in\R^+$
%
\begin{eqnarray}
\label{6_bis} \hspace*{-30pt}\Es_\by\bigl[H^{\sigma_{k-1}}(T \wedge
\tau_1\wedge\sigma_k)+H^{\sigma
_{k-1}}(
\tau_1\wedge\sigma_k) |\Ff_{\sigma_{k-1}} \bigr]
\le2H^{\sigma_{k-1}}(\sigma_{k-1})
\nonumber
\\[-8pt]
\\[-8pt]
\eqntext{\mbox{on the set } \{T \wedge
\tau_1>\sigma_{k-1}\}.}
\end{eqnarray}
\end{proposition}

Once this proposition is proven, by (\ref{5}) we will have
\begin{eqnarray*}
\Es_\by\bigl( H(\tau_1\wedge\sigma_k)
\bigr) &\le&\Es_\by\bigl( H(\tau_1\wedge
\sigma_{k-1})\bigr),
\\
\Es_\by\bigl( H(T \wedge\tau_1\wedge
\sigma_k) + H(\tau_1\wedge\sigma_k)
\bigr) &\le&\Es_\by \bigl(H(T \wedge\tau_1\wedge
\sigma_{k-1}) + H(\tau_1\wedge\sigma_{k-1})
\bigr)
\end{eqnarray*}
and iterating these inequalities we obtain Proposition~\ref
{prop_energy_clusters} because $H(0)=V(\by)$.
In order to prove Proposition~\ref{lem_recur}, we have to consider two cases.

\begin{pf*}{Proof of Proposition~\ref{lem_recur} when $\bolds{k}$ is odd
(i.e., $\bolds{|U_{23}|^2}$ goes downhill)}
Suppose $\tau_1> \sigma_{k-1}$ and look at the dynamics during the
interval $[\sigma_{k-1}, \sigma_k\wedge\tau_1)$. This case is simple
because no balls can
collide, hence the local time term $L_{23}$ in (\ref{4}) vanishes.
Therefore by (\ref{4}), we have, on this time interval,
\begin{eqnarray*}
\mathrm{d}H^{\sigma_{k-1}}(t)&=&6\mathrm{e}^{\lambda(t-\sigma_{k-1})}
\bigl(U_1(t), \mathrm{d}B_1(t)\bigr)+6\mathrm{e}^{\lambda(t-\sigma_{k-1})}
\bigl(U_{23}(t), \mathrm{d}B_{23}(t)\bigr) \qquad\mbox{(martingale part)}
\\
&&{}+ H^{\sigma_{k-1}}(t) \biggl( \frac{6d}{V(Y(t))} +\lambda-6a
\biggr) \,\mathrm{d}t.
\end{eqnarray*}
On $[\sigma_{k-1}, \sigma_k\wedge\tau_1)$ the border level $R(t)$
equals $R$, so we have
$V(Y(t))>2(R+1)+2 (\frac{1}{2} +\medm) = 2R+3+2\medm$. Therefore,
\begin{eqnarray*}
\Es_\by\bigl[H^{\sigma_{k-1}}(
\tau_1\wedge\sigma_k) |\Ff_{\sigma
_{k-1}} \bigr]&\leq&
H^{\sigma_{k-1}}(\sigma_{k-1})
\\
&&{}+ \biggl(\frac{6d}{2R+3+2\medm}+\lambda-6a \biggr)\Es_\by\biggl
[\int
_{\sigma
_{k-1}}^{\tau_1\wedge\sigma_k}H^{\sigma_{k-1}}(t) \,\mathrm{d}t \Big|\Ff
_{\sigma
_{k-1}} \biggr].
\end{eqnarray*}
Since $ H^{\sigma_{k-1}}(\sigma_{k-1})=V(Y(\sigma_{k-1}))$, this yields
(\ref{6}) provided that
%
\begin{equation}
\label{cond_lambda} \lambda\le6a-\frac{6d}{2R+3+2\medm}.
\end{equation}
Note that for any fixed time horizon $T\in\R^+$, the above calculation
also holds with $\tau_1$ replaced by $\tau_1 \wedge T$.
\end{pf*}

\begin{pf*}{Proof of Proposition~\ref{lem_recur} when $\bolds{k}$ is even
(i.e., $\bolds{|U_{23}|^2}$ goes uphill)}
We look at the dynamics during the interval $[\sigma_{k-1}, \sigma
_k\wedge\tau_1)$ again.

Up to a martingale term, $H^{\sigma_{k-1}}(\tau_1\wedge\sigma_k)
-H^{\sigma_{k-1}}(\sigma_{k-1})$ is equal to
\begin{eqnarray*}
&&(\lambda-6a)\int_{\sigma_{k-1}}^{\tau_1\wedge\sigma_k} 3
\mathrm{e}^{\lambda(
s-\sigma_{k-1})}
\bigl|U_1(s)\bigr|^2 \,\mathrm{d}s +\frac{3d}{\lambda} \bigl(
\mathrm{e}^{\lambda( \tau_1\wedge\sigma
_k-\sigma
_{k-1})} -1 \bigr) \\
&&\qquad{}+3 \mathrm{e}^{\lambda( \tau_1\wedge\sigma_k-\sigma
_{k-1})}\bigl|U_{23}(\tau
_1\wedge\sigma_k)\bigr|^2 -3\bigl|U_{23}(
\sigma_{k-1})\bigr|^2.
\end{eqnarray*}
In this case, on $[\sigma_{k-1}, \sigma_k\wedge\tau_1)$, thanks to
(\ref{eq:boundsmedians}), $|U_1(s)|^2 \ge
\frac{1}{3}(1+R')-\frac{1}{3}(\frac{1}{2}+\medp)$.

Moreover, $|U_{23}(\sigma_{k-1})|^2=\frac{1}{2}+\medm$ and
$|U_{23}(\tau
_1\wedge\sigma_k)|^2 \le\frac{1}{2}+\medp$.

Thus for any $\lambda<6a$
%
\begin{eqnarray}
\label{Hdecreases} &&\Es_\by\bigl( H^{\sigma_{k-1}}(\tau_1
\wedge\sigma_k) - H^{\sigma
_{k-1}}(\sigma_{k-1}) |
\Ff_{\sigma_{k-1}} \bigr)
\nonumber
\\
&&\quad \le \biggl( (\lambda-6a) \biggl( \frac{1}{2}+R'-\medp
\biggr) +3d \biggr) \Es_\by\biggl( \frac{ \mathrm{e}^{\lambda( \tau
_1\wedge\sigma_k-\sigma_{k-1})}
-1 }{\lambda} \Big|
\Ff_{\sigma_{k-1}} \biggr)
\nonumber
\\[-8pt]
\\[-8pt]
\nonumber
&&\qquad{}+ 3\biggl(\frac{1}{2}+\medp\biggr) \Es_\by\bigl(
\mathrm{e}^{\lambda( \tau_1\wedge
\sigma
_k-\sigma_{k-1})} | \Ff_{\sigma_{k-1}} \bigr) - 3\biggl(\frac
{1}{2}+\medm
\biggr)
\\
&&\quad= \bigl( R'(\lambda-6a) + 2\lambda(1+\medp) -3a(1-2\medp) +3d
\bigr) \Es_\by\biggl( \frac{ \mathrm{e}^{\lambda( \tau_1\wedge\sigma
_k-\sigma_{k-1})}
-1 }{\lambda} \Big| \Ff_{\sigma_{k-1}} \biggr)
+3(\medp-\medm).\nonumber
\end{eqnarray}
The key point in the whole proof is the fact that, under an appropriate
choice of the parameters,
this last expectation is finite and admits a uniform lower
bound.

\begin{lemma}\label{minor_even} 
There exists $\underline{C}$ depending only on $\medp$, $\medm$ such
that for each even $k$ and for $\lambda$ small enough
\[
0 < \underline{C} \le\Es_\by\biggl( \frac{ \mathrm{e}^{\lambda( \tau
_1\wedge
\sigma_k-\sigma_{k-1})} -1 }{\lambda} \Big|
\Ff_{\sigma_{k-1}} \biggr) <+\infty.
\]
\end{lemma}

%
Once this lemma is proved, there will be an $R'$ large enough (hence an
$R$ large enough) for
%
\begin{equation}
\label{cond_R'} \bigl(R'(\lambda-6a) +2\lambda(1+\medp) -3a(1-2
\medp) +3d\bigr)\underline{C}+3(\medp-\medm) \le0
\end{equation}
to hold and this will imply for $\tau_1>\sigma_{k-1}$
%
\begin{eqnarray}
\label{Hdecreases_bis} \Es_\by\bigl( H^{\sigma_{k-1}}(\tau_1
\wedge\sigma_k) - H^{\sigma
_{k-1}}(\sigma_{k-1}) |
\Ff_{\sigma_{k-1}} \bigr) \le0,
\end{eqnarray}
which rewrites into (\ref{6}).

Note that, as in the previous case, calculation (\ref{Hdecreases}) also
holds with $\tau_1$ replaced by $\tau_1 \wedge T$ for any fixed time
horizon $T\in\R^+$.
Summing the expressions with and without finite time horizon, and using
the lower bound $0$ for
$\Es_\by( \frac{ \mathrm{e}^{\lambda( \tau_1\wedge\sigma_k\wedge T
-\sigma
_{k-1})} -1 }{\lambda} | \Ff_{\sigma_{k-1}} )$, we obtain
that on
$T\wedge\tau_1>\sigma_{k-1}$
\[
\Es_\by\bigl( H^{\sigma_{k-1}}(\tau_1\wedge
\sigma_k)+H^{\sigma
_{k-1}}(T\wedge\tau_1\wedge
\sigma_k) -2H^{\sigma_{k-1}}(\sigma_{k-1}) |
\Ff_{\sigma_{k-1}} \bigr) \le0
\]
as soon as
%
\begin{equation}
\label{cond_R'_bis} \bigl(R'(\lambda-6a) +2\lambda(1+\medp) -3a(1-2
\medp) +3d\bigr)\underline{C} +6(\medp-\medm) \le0.
\end{equation}
\end{pf*}

Clearly, we can forget about condition (\ref{cond_R'}) as any set of
parameter satisfying condition~(\ref{cond_R'_bis}) will satisfy (\ref
{cond_R'}) too.
From now on, our aim is to prove Lemma~\ref{minor_even}. The finiteness
of $\Es_\by( \frac{ \mathrm{e}^{\lambda( \tau_1\wedge\sigma
_k-\sigma
_{k-1})} -1 }{\lambda}
| \Ff_{\sigma_{k-1}} )$ is obtained in Section~\ref
{finite_exp_moment} and the uniform lower bound is constructed in
Section~\ref{subsec_lower_bound}.

\subsection{Existence of exponential moments of \texorpdfstring{$\tau_1\wedge\sigma_k-\sigma_{k-1}$}{the stopping time} 
for even~\texorpdfstring{$k$}{k}}\label{finite_exp_moment}

We need a proof that for small enough $\lambda$'s the exponential moment
$\Es_\by( \mathrm{e}^{\lambda( \tau_1\wedge\sigma_k-\sigma_{k-1})}
|
\Ff_{\sigma_{k-1}} )$ is finite.
We use a comparison argument. Since a similar comparison argument will
be needed to obtain a lower bound on the exponential moment, we
directly construct a
double inequality, though an upper bound is sufficient for our purpose
in this section.

\subsubsection{Comparison with the level hitting time of a simple
reflected SDE}
Consider a Wiener process $B^U$, independent on $W_1$, such that
\[
B^U(t)=\int_0^t
\frac{1}{|U_{23}(s)|} \bigl(U_{23}(s), \mathrm{d}B_{23}(s)\bigr),\qquad t<
\tau_1 ,
\]
and the process $U$ solution to the following one-dimensional SDE with
reflection at the point $\frac{1}{2}$:
\[
U(t)=\bigl|U_{23}(0)\bigr|^2 +\int_0^t
\bigl(d-6aU(s)\bigr) \,\mathrm{d}s +2\int_0^t
U^{1/2}(s) \,\mathrm{d}B^U(s) +L^U(t).
\]
Then the processes $|U_{23}(t)|^2$ and $U(t)$ coincide up to the time
$\tau_1$, and $L^U(t)=2L_{23}(t)$ for $t<\tau_1$.
It is sufficient to prove the finiteness of $\Es_{\frac{1}{2}+\medm} (
\mathrm{e}^{\lambda\sigma})$ where $\sigma=\inf\{t\dvt U(t)=\frac
{1}{2}+\medp\}$.
We make a time change, that is, we put
\[
\zeta_t=\int_0^t 4 U(s)\,\mathrm{d}s,\qquad
\chi_t=\inf\{r\dvt\zeta_r\geq t\},\qquad \tilde U(t)=U(
\chi_t).
\]
Then $\tilde U$ satisfies the one-dimensional SDE with reflection at
the point $\frac{1}{2}$
\[
\mathrm{d}\tilde U(t)=\frac{d -6a \tilde U(t)}{4 \tilde U(t)} \,\mathrm{d}t
+ \mathrm{d}\tilde
B(t)+\mathrm{d}\tilde L(t),
\]
where $\tilde B$ is a Wiener process. Since $\frac{1}{2} \le\tilde
U(t) \le\frac{1}{2}+\medp$ up to time $\sigma$, then
\[
(2+4\medp)\sigma\ge\tilde\sigma:=\zeta_\sigma= \inf\bigl\{
t\dvt\tilde
U(t)=\tfrac{1}{2}+\medp\bigr\} \ge\sigma.
\]
Since the drift of $\tilde U(t)$ is bounded from above and from below
by some constants
\[
C_1:= -\frac{3}{2}a \le\frac{d -6a \tilde U(t)}{4 \tilde U(t)} \le
\frac{d}{2}-\frac{3}{2}a < \frac{d}{2} =: C_2 ,
\]
we can compare $\tilde U(t) $ with reflected Brownian motions with
constant drifts $C_1$ and $C_2$, as in Ward and Glynn \cite{WG},
Proposition~2.
We then obtain $\hat U_1 \le\tilde U \le\hat U_2$ where $\hat U_i,
i=1,2$ satisfy the one-dimensional SDEs with reflection at the point
$\frac{1}{2}$
\[
\mathrm{d}\hat U_i(t)=C_i \,\mathrm{d}t+ \mathrm{d}\hat B(t)+\mathrm{d} \hat
L_i(t) ,\qquad \hat U_i(0)=\tfrac
{1}{2}+\medm,
\]
where $\hat B$ is an $\R$-valued Brownian motion and $\hat L_i(t)=\int
_0^t \un_{\hat U_i(s)=\frac{1}{2}}\,\mathrm{d}\hat L_i(s)$.
This allows us to compare the $\delta$-level hitting times of the three
processes:
\[
\hat\sigma_1:=\inf\bigl\{t\dvt\hat U_1(t)\ge
\tfrac{1}{2}+\medp\bigr\} \ge\tilde\sigma\ge\hat\sigma_2:=
\inf\bigl\{t\dvt\hat U_2(t)\ge\tfrac{1}{2}+\medp\bigr\}
\]
and we obtain
%
\begin{equation}
\label{compar_hitting_times} \frac{\hat\sigma_2}{2+4\medp} \le
\sigma\le\hat\sigma_1.
\end{equation}
In the sequel, we compute the exponential moments of hitting times
$\hat\sigma_i, i\in\{1,2\}$. For the time being, we drop the indices
on $\hat\sigma_i$,
$C_i$ and $\hat L_i$.

\subsubsection{Exponential moments of level hitting times}

It is equivalent to consider the hitting time of $\frac{1}{2}+\medp$
for a Brownian motion starting from $\frac{1}{2}+\medm$ with constant
drift $C$
and reflection at $1/2$ or to consider the hitting time of $\medp$ for
a Brownian motion starting from $\medm$ with constant drift $C$ and
reflection at $0$.
Girsanov theorem for processes with reflection Kinkladze \cite
{Kinkladze82} and
Doob's optional sampling theorem implies that for all negative~$\lambda$
\[
\Es_{{1}/{2}+\medm} \bigl( \mathrm{e}^{\lambda\hat\sigma} \bigr) =
\mathrm{e}^{C(\medp-\medm)} \Es\bigl(
\mathrm{e}^{(\lambda-{C^2}/{2})\inf\{
t;|\medm
+\hat B(t)|=\medp\}} \mathrm{e}^{-({C}/{2})\lim_{\eps\to0}({1}/{\eps
})\int_0^t\un
_{[0;\eps
[}(|\medm+\hat B(t)|)} \bigr).
\]
Suppose $C\neq0$.
Then using Formula 2.3.3 in Borodin and Salminen \cite{BoSa} for $r=0,
x=\medm, z=\medp,
\alpha= C^2/2 - \lambda, \gamma= C/2$, which holds for any $\lambda<C^2/2$,
%
\begin{equation}
\label{def_psi} \Es_{{1}/{2}+\medm} \bigl( \mathrm{e}^{\lambda\hat
\sigma} \bigr) =
\mathrm{e}^{C(\medp-\medm)} \frac{v \cosh(C \medm v) +\sinh(C \medm v)}{v
\cosh
(C \medp v) +\sinh(C \medp v)} =: \Psi(v),
\end{equation}
where $v(\lambda):=\sqrt{1-2\frac{\lambda}{C^2}}$.
This is an analytical function of $\lambda$ thus the formula holds
as long as $v(\lambda)$ is well defined and the
denominator does not vanish. But any positive $v$ such that the
denominator vanishes satisfies $\frac{x\cosh(x)}{\sinh(x)}=-C \medp$
for $x=C \medp v$. Since the
function $x\mapsto\frac{x\cosh(x)}{\sinh(x)}$ is larger than $1$
on the whole $\R$, the condition $-C \medp< 1$ ensures that the
$\lambda$-exponential
moment exists for positive $\lambda$'s satisfying $\lambda<C^2/2$.

$\sigma\le\tilde\sigma\le\hat\sigma_1$ thus $\Es_{{1}/{2}+\medm} ( \mathrm{e}^{\lambda\sigma} )$ is finite as soon as
$\lambda<C_1^2/2$ and $-C_1
\medp< 1$, hence
%
\begin{equation}
\label{cond_lambda_bis} \lambda<\frac{9}{8}a^2 \quad\mbox{and}\quad \medp<
\frac{2}{3a} \quad\Longrightarrow\quad\Es_\by\bigl( \mathrm{e}^{\lambda( \tau
_1\wedge\sigma_k-\sigma_{k-1})} |
\Ff_{\sigma_{k-1}} \bigr) <+\infty.
\end{equation}
From now on, we assume $\lambda<\frac{9}{8}a^2$ and $\medp< \frac{2}{3a}$.

\subsection{Lower bound for the exponential moment of \texorpdfstring
{$\tau_1\wedge\sigma_k-\sigma_{k-1}$}{the stopping time} for even
\texorpdfstring{$k$}{k}}\label{subsec_lower_bound}
We replace in the definition of $\Es_\by( \frac{ \mathrm{e}^{\lambda(
\tau
_1\wedge\sigma_k-\sigma_{k-1})} -1 }{\lambda} | \Ff_{\sigma_{k-1}}
)$
the stopping time $\tau_1$, which is expressed in the terms of the
minimum of $|Y_1-Y_2|^2$ and $ |Y_1-Y_3|^2$, by another one, expressed
in the terms of
$U_1$.

Note that if $\tau_1 \in[\sigma_{k-1}, \sigma_k)$ with $k$ even, then
$|U_{23}(\tau_1)|^2 \le\frac{1}{2}+\medp$ and thanks to (\ref
{eq:boundsmedians})
\[
\bigl|U_1(\tau_1)\bigr|^2 \le\frac{4}{3}
\bigl(1+R'\bigr)+\frac{2}{3}\biggl(\frac{1}{2}+\medp
\biggr)= \frac{5+4R'+2\medp}{3}.
\]
Thus, $\tau_1\wedge\sigma_k\geq\rho_k\wedge\sigma_k$ for $\rho
_k=\inf
\{t\geq\sigma_{k-1}; |U_1(t)|^2\le\frac{5+4R'+2\medp}{3}\}$.
Observe that, because we have assumed that $\tau_1>\sigma_{k-1}$,
equality (\ref{eq:boundsmedians}) also implies
\[
\bigl|U_1(\sigma_{k-1})\bigr|^2 \ge
\tfrac{2}{3}(1+R) -\tfrac{1}{3}\bigl( \tfrac{1}{2} +\medm
\bigr)= \tfrac{1}{3} \bigl( \tfrac{3}{2} +2R-\medm\bigr).
\]
Using the fact that $(\mathrm{e}^{\lambda s}-1)/\lambda\ge s$ for $s,\lambda
>0$, we have
%
\begin{eqnarray}
\label{minor_esp} \Es_\by\biggl( \frac{ \mathrm{e}^{\lambda( \tau_1\wedge
\sigma_k-\sigma_{k-1})}
-1 }{\lambda} \Big|
\Ff_{\sigma_{k-1}} \biggr) &\ge& \Es_\by(\tau_1\wedge
\sigma_k-\sigma_{k-1} | \Ff_{\sigma
_{k-1}} )
\nonumber
\\[-8pt]
\\[-8pt]
\nonumber
&\ge& \inf\Es\biggl( \rho\wedge\sigma| U_{1}(0)=u_{1},
|U_{23}|^2(0)=\frac{1}{2}+\medm\biggr),
\label{21}
\end{eqnarray}
where the infimum is taken among all initial conditions $u_1$ such that
$|u_1|^2\geq\frac{1}{3}(\frac{3}{2} +2R-\medm)$,
\[
\rho=\inf\biggl\{t\geq0\dvt\bigl|U_1(t)\bigr|^2 \le
\frac{5+4R'+2\medp}{3}\biggr\}
\]
and
\[
\sigma=\inf\bigl\{t\geq0\dvt\bigl|U_{23}(t)\bigr|^2 \ge
\tfrac{1}{2}+\medp\bigr\}.
\]
Because $U_1$ and $U_{23}$ are independent up to time $\tau_1$, we can
estimate the right-hand side in (\ref{21}) in the following way:
for an arbitrary $Q>0$ which will be chosen later,
%
\begin{eqnarray}
\label{22} &&\inf_{|u_1|^2\ge({1}/{3})({3}/{2} +2R-\medm)}\Es
\biggl[ \rho\wedge\sigma|
U_{1}(0)=u_{1}, |U_{23}|^2(0)=
\frac{1}{2}+\medm\biggr]
\nonumber
\\[-8pt]
\\[-8pt]
\nonumber
&&\qquad\ge\Es\biggl(\sigma\wedge Q | |U_{23}|^2(0)=
\frac{1}{2}+\medm\biggr) \inf_{|u_1|^2\ge({1}/{3})({3}/{2}
+2R-\medm)} \Ps\bigl(\rho>Q |
U_{1}(0)=u_{1}\bigr).
\nonumber
\end{eqnarray}
Let us compute a lower bound for each factor.

\subsubsection{Lower bound for \texorpdfstring{$\Ps(\rho>Q)$}{$P(rho>Q)$}}
By It\^o formula $ \mathrm{d}|U_1(t)|^2=2(U_1(t), \mathrm{d} B_1(t))-6a |U_1(t)|^2
\,\mathrm{d}t+d \,\mathrm{d}t$,
and $ \mathrm{d}\log(|U_1(t)|^2)=2|U_1(t)|^{-2}(U_1(t), \mathrm{d} B_1(t))-6a
\,\mathrm{d}t+(d-2)|U_1(t)|^{-2} \,\mathrm{d}t$.

Denote
\[
M_t=2\int_0^t\bigl|U_1(s)\bigr|^{-2}
\bigl(U_1(s), \mathrm{d} B_1(s)\bigr),
\]
then, for $d \ge2$,
\[
\log\bigl(\bigl|U_1(t)\bigr|^2\bigr)\ge\log\bigl(\bigl|U_1(0)\bigr|^2
\bigr)+M_t-6at.
\]
Note that $|U_1(s)|^2\ge\frac{5+4R'+2\medp}{3}$ up to time $\rho$ thus
\[
\Es\bigl( M_{t \wedge\rho}^2 \bigr) = 4\Es\int_0^{t \wedge\rho}
\bigl|U_1(s)\bigr|^{-2} \,\mathrm{d}s\le\frac{12 t}{5+4R'+2\medp},
\]
so that, by the Doob inequality,
\[
\Ps\biggl( \sup_{s\le Q}|M_{s \wedge\rho}| \ge\sqrt{
\frac
{24 Q}{5+4R'+2\medp}} \biggr) \le\frac{5+4R'+2\medp}{24 Q} \Es
\bigl( M^2_{Q \wedge\rho}
\bigr) \le\frac{1}{2}.
\]
Then, with probability at least $1/2$,
\begin{eqnarray*}
\inf_{s\le Q} \log\bigl(\bigl|U_1(s \wedge
\rho)\bigr|^2\bigr)&\ge&\log\bigl(\bigl|U_1(0)\bigr|^2\bigr)-
\sqrt{\frac{24 Q}{5+4R'+2\medp}}-6aQ\\
& \ge&\log\bigl(\bigl|U_1(0)\bigr|^2
\bigr)-\sqrt{6Q/R'}-6aQ.
\end{eqnarray*}
This means that
%
\begin{equation}
\label{24} \Ps\bigl(\rho>Q| U_{1}(0)=u_{1}\bigr)
\ge1/2\qquad \forall|u_1|^2\ge\tfrac
{1}{3}\bigl(
\tfrac{3}{2}+2R-\medm\bigr)
\end{equation}
holds true as soon as (large) $R$, $Q$ and $R'$ are chosen in such a
way that
\[
\log\biggl( \frac{1}{3}\biggl(\frac{3}{2}+2R-\medm\biggr) \biggr)-
\sqrt{6Q/R'}-6aQ > \log\biggl(\frac
{5+4R'+2\medp}{3} \biggr)
\]
that is,
%
\begin{equation}
\label{cond_R} \tfrac{3}{2}+2R-\medm> \mathrm{e}^{\sqrt{6Q/R'}+6aQ}\bigl(
5+4R'+2\medp\bigr).
\end{equation}

\subsubsection{Lower bound for \texorpdfstring{$\Es(\sigma\wedge Q)$}{$E(sigma wedge Q)$}}
We have $\sigma\wedge Q = \sigma-(\sigma-Q)\un_{\sigma> Q} \ge
\sigma
-(\sigma-Q)\frac{\sigma}{Q}$ hence
\[
\Es(\sigma\wedge Q) \ge2\Es(\sigma)-\frac{1}{Q} \Es\bigl(
\sigma^2\bigr).
\]
The comparison argument developed in Section~\ref{finite_exp_moment}
leads to
%
\begin{equation}
\label{minor_sigma_Q} \Es(\sigma\wedge Q) \ge\frac{2}{2+4\medp}
\Es(\hat
\sigma_2) -\frac
{1}{Q} \Es\bigl(\hat\sigma_1^2
\bigr).
\end{equation}
We need a lower bound for the first moment of $\hat\sigma_2$ and an
upper bound for the second moment of~$\hat\sigma_1$.
To this end, we use the exponential moment of $\hat\sigma$ given by
(\ref{def_psi}) for $C\neq0$ with $-C\medp<1$, on a neighbourhood of
zero for $\lambda$.
Differentiating twice in (\ref{def_psi}) at $\lambda=0$, we obtain the
first and second moment of $\hat\sigma$.
In order to simplify the derivative computations, from now on we make
the simplifying choice
\[
\medm=\medp/2.
\]
We obtain
\[
\Es_{{1}/{2}+\medm} ( \hat\sigma) = -\frac{\Psi'(1)}{C^2}
=\frac{\mathrm{e}^{-2C\medp}-\mathrm{e}^{-C\medp}}{2C^2}+
\frac{\medp}{2C} 
= \frac{3}{4}\medp^2 -
\frac{7}{12}C\medp^3 + \frac{1}{2C^2} \sum
_{k=4}^{+\infty} \frac{(-C\medp)^k}{k!}\bigl(2^k-1
\bigr).
\]
The right-hand side series is positive as soon as $C \medp\le2$
because the
sequence $u_k:=\frac{2^k-1}{k!}$ satisfies $u_k > 2 u_{k+1}$ for all
$k\ge4$.
So, since $C_2=\frac{d}{2}$, if $\medp\le\frac{4}{d}$
\[
\Es_{{1}/{2}+\medm} ( \hat\sigma_2 ) \ge\frac
{3}{4}
\medp^2 -\frac{7}{12}C_2\medp^3 =
\frac{\medp^2}{12}\biggl(9-7\frac
{d}{2}\medp\biggr)
\]
and in particular
\[
\Es_{{1}/{2}+\medm} ( \hat\sigma_2 ) \ge\frac
{\medp
^2}{6} \quad\mbox{as
soon as} \quad\medp\le\frac{2}{d}.
\]
Moreover,
\[
\Es_{{1}/{2}+\medm} \bigl( \hat\sigma^2 \bigr) = \frac
{\Psi
''(1)-\Psi'(1)}{C^4}
=\frac{-1}{C^3} \int_\medm^\medp
\bigl(\mathrm{e}^{-2Cx}-1\bigr) \bigl(\mathrm{e}^{-2C\medp}+2C\medp+1\bigr)+2Cx
\bigl(\mathrm{e}^{-2Cx}+1\bigr) \,\,\mathrm{d}x.
\]
For any negative $C_1$ such that $-C_1 \medp< 1$, one has
\[
\Es_{{1}/{2}+\medm} \bigl( \hat\sigma_1^2 \bigr) \le
\frac
{-1}{C_1^3} \int_\medm^\medp
\bigl(\mathrm{e}^{-2C_1x}-1\bigr) \bigl(\mathrm{e}^{-2C_1\medp}+1\bigr) \,\,\mathrm{d}x \le
\frac{\medp}{2(-C_1)^3} \bigl(\mathrm{e}^{-4C_1\medp}-1\bigr) \le2\frac
{\medp
^2}{C_1^2}\mathrm{e}^4
\]
because $\mathrm{e}^{4x}-1 \le4\mathrm{e}^4 x$ for $x$ between $0$ and $1$.

Since $C_1=\frac{-3a}{2}$, inequality (\ref{minor_sigma_Q}) leads to
%
\begin{equation}
\label{minor_esp_sigma} \Es_{{1}/{2}+\medm}(\sigma\wedge Q)
\ge\frac{\medp
^2}{6(1+2\medp
)}-
\frac{8\medp^2}{9a^2 Q}\mathrm{e}^4
\end{equation}
under the conditions that $\medp\le\min(\frac{2}{3a},\frac{2}{d})$.

Using (\ref{minor_esp}), (\ref{22}), (\ref{24}) and (\ref
{minor_esp_sigma}), we have obtained
%
\begin{equation}
\label{3.25_bis} \Es_\by\biggl( \frac{ \mathrm{e}^{\lambda( \tau_1\wedge
\sigma_k-\sigma_{k-1})}
-1 }{\lambda}\Big |
\Ff_{\sigma_{k-1}} \biggr) \ge\frac{\medp^2}{2} \biggl( \frac
{1}{6+12\medp} -
\frac
{8\mathrm{e}^4}{9a^2 Q} \biggr).
\end{equation}
This induces our choice of $Q$ to simplify the right-hand side of
(\ref{3.25_bis}):
\[
Q=\frac{32(1+2\medp)\mathrm{e}^4}{3a^2} \quad\Leftrightarrow\quad\frac
{8\mathrm{e}^4}{9a^2 Q} = \frac{1}{12+24\medp}
\]
and we obtain the lower bound
%
\begin{equation}
\label{cond_Delta} \Es_\by\biggl( \frac{ \mathrm{e}^{\lambda( \tau
_1\wedge\sigma_k-\sigma_{k-1})}
-1 }{\lambda} \Big|
\Ff_{\sigma_{k-1}} \biggr) \ge\frac{\medp
^2}{24(1+2\medp)} \qquad\mbox{if } \medp<
\frac{2}{3a} \mbox{ and } \medp\le\frac{2}{d}.
\end{equation}
%
\subsection{Choice of the parameters} 
Recall that $\medm=\medp/2$.
We have to choose four parameters $\medp, R, R'$ and $\lambda$, which
should satisfy the following five
conditions:
\begin{eqnarray*}
&&\medp< \frac{2}{3a} \quad\mbox{and}\quad \medp\le\frac{2}{d} \qquad\mbox{from
(\ref{cond_Delta})},
\\
&&\lambda< \frac{9}{8}a^2 \qquad\mbox{from (\ref{cond_lambda_bis})},
\\
&&\lambda\le6a-\frac{6d}{2R+3+\medp} \qquad\mbox{from (\ref{cond_lambda})},
\\
&&\bigl( R'(\lambda-6a) +2\lambda(1+\medp) -3a(1-2\medp) +3d \bigr)
\frac
{\medp^2}{24(1+2\medp)} +3\medp\le0
\end{eqnarray*}
that is,
\begin{eqnarray*}
&&R'(6a-\lambda) -2\lambda(1+\medp) +3a(1-2\medp)
-3d \ge72\biggl(\frac{1}{\medp}+2\biggr) \qquad\mbox{from (\ref{cond_R'_bis})},
\\
&&\frac{3-\medp}{2}+2R > \mathrm{e}^{\sqrt{6Q/R'}+6aQ}\bigl( 5+4R'+2\medp
\bigr)\qquad
\mbox{for } Q=\frac{32(1+2\medp)\mathrm{e}^4}{3a^2} \qquad\mbox{from (\ref{cond_R})}.
\end{eqnarray*}
We choose $\medp=\frac{2}{3a+d} \le1$ which complies with (\ref
{cond_Delta}). We fix $\lambda=\min(a,a^2)$.
Condition (\ref{cond_R'_bis}) is satisfied as soon as $R'(6a-\lambda)
\ge72(\frac{3a+d}{2}+2)+3d-3a+6a\medp+2\lambda+\frac{4\lambda}{3a+d}$.
Since $6a\medp\le4$ and $\lambda\le a$, (\ref{cond_R'_bis}) holds in
particular if
\[
R'=\frac{22a+8d+30}{a}.
\]
The last parameter $R$ will be taken large enough to satisfy $2R \ge
\mathrm{e}^{\sqrt{6Q/R'}+6aQ}(7+4R')$ which implies (\ref{cond_R}).
First, remark that $R'>22$ with our choice, hence
\[
\sqrt{6Q/R'}+6aQ \le\sqrt{\frac{64(1+2\medp)\mathrm{e}^4}{22a^2}}+\frac
{64(1+2\medp)\mathrm{e}^4}{a}
\le\frac{10\,505}{a}.
\]
Noticing that $7+4R' \le(95a+32d+120)/a$ with the choice of $R'$ we
made, we obtain a sufficient condition for (\ref{cond_R}) to hold:
\[
R \ge\frac{48a+16d+60}{a} \mathrm{e}^{10\,505/a}.
\]
Such an $R$ satisfies $R>16d/a$ hence is more than sufficient for (\ref
{cond_lambda}) to hold.

This completes the proofs of Lemma~\ref{minor_even} and Proposition~\ref
{lem_recur}, hence Proposition~\ref{prop_energy_clusters} holds.
This in turn completes
the proof of Proposition~\ref{prop_sect_3}.

\begin{appendix}\label{app}
\section*{Appendix: $D$ has a Lipschitz
boundary}\label{appendix_Lipschitz}
The following lemma is useful to apply results from Bass and Hsu \cite
{BassHsu}, Chen, Fitzsimmons and Williams \cite{CFW}, Fukushima and
Tomisaki \cite{Fuktom} to the hard ball process $X$.

%
\begin{lemmaa}\label{lemLipschitz}
The domain
\[
D = \bigl\{\bx\in\bigl(\mathbb R^d\bigr)^n ;
|x_i - x_j|>r \mbox{ for all } i \neq j \bigr\}
\]
has a Lipschitz boundary.
\end{lemmaa}
%

\begin{pf}
Define the function $f_{ij}$ on $(\mathbb R^d)^n$ by $f_{ij}(\bx)=|x_i
- x_j|^2-r^2$.
Fix $\bx\in\partial D$. Proceeding like in the proof of Proposition~4.1
in Fradon \cite{FradonGlob}, one can show that
there exits a unit vector $\bv$ such that, for each pair ($i,j$) of
colliding balls of $\bx$, $\nabla
f_{ij}(\bx)\cdot\bv\ge\frac{r}{n\sqrt{2n}}>0$. Indeed $\bv$ is the
direction in which each
colliding ball goes away from the gravity center of the collision.

Take $m:=n d$. By continuity,
\[
\eps(\bx)=\inf\bigl\{ \bigl|\bx'-\bx\bigr| , \bx'\in\partial D
\mbox{ and } \exists(i,j) \mbox{ s.t. } |x_i - x_j|>r
\mbox{ and }\bigl |x'_i - x'_j\bigr|=r
\bigr\}
\]
is positive. On the ball with center $\bx$ and radius $\eps(\bx)$, we
choose an orthonormal coordinate system $(y_1,\ldots,y_m)$ with point
$\bx$
as the origin and direction $\bv$ as the last axis, that is, $\bx'$ has
coordinates $(y_1,\ldots,y_m)$ with $y_m=(\bx'-\bx)\cdot\bv$.

Let us write $f_{ij} \circ h$ for the function $f_{ij}$ expressed in
this coordinates system. The partial derivative of $f_{ij}$ at the
origin with respect to
the last coordinate is given by
\[
\frac{\partial(f_{ij} \circ h)}{\partial y_m}(0) =\lim_{\eta\to0}
\frac{f_{ij}(\bx+\eta\bv)-f_{ij}(\bx)}{\eta}=\nabla
f_{ij}(\bx)\cdot\bv>0.
\]
Therefore, due to the implicit function theorem, there exists a $\CC
^1$-function $g_{ij}$ on $\R^m$ such that
$f_{ij}(h(y_1,\ldots,g_{ij}(y_1,\ldots,y_{m-1},\cdot)))=\mathrm{Id}$ for each
$h(y_1,\ldots,y_m) \in B(\bx,\eps)$ where $\eps<\eps(\bx)$ is such
that, on
$B(\bx,\eps)$, all the functions $\nabla f_{ij}(\cdot).\bv$ stay
positive for any pair ($i,j$) of colliding balls of $\bx$.
The maps $y_m \mapsto f_{ij}(h(y_1,\ldots,y_{m}))$ are increasing,
so that for $h(y_1,\ldots,y_m)$ in $B(\bx,\eps)$:
\begin{eqnarray*}&& y_m > \max\bigl\{
g_{ij}(y_1,\ldots,y_{m-1},0) \mbox{ s.t. }
|x_i - x_j|=r \bigr\}
\\
&&\quad\Leftrightarrow\quad\forall i<j \quad\mbox{s.t.}\quad |x_i -
x_j|=r, y_m > g_{ij}(y_1,\ldots
,y_{m-1},0)
\\
&&\quad\Leftrightarrow\quad\forall i<j\quad\mbox{s.t.}\quad |x_i -
x_j|=r,\\
&&\hspace*{14pt}\qquad f_{ij}\bigl(h(y_1,\ldots
,y_{m-1},y_m)\bigr)>f_{ij}\bigl(h
\bigl(y_1,\ldots,g_{ij}(y_1,
\ldots,y_{m-1},0)\bigr)\bigr)=0
\\
&&\quad\Leftrightarrow\quad\forall i<j\ f_{ij}\bigl(h(y_1,
\ldots,y_{m-1},y_m)\bigr)>0
\\
&&\quad\Leftrightarrow\quad h(y_1,\ldots,y_m) \in D.
\end{eqnarray*}
Note that, since the $g_{ij}$ are $\CC^1$, they are Lipschitz
continuous with Lipschitz constant $C_{ij}$,
which leads to the Lipschitz continuity of the function
\[
(y_1,\ldots,y_{m-1}) \mapsto\max\bigl\{
g_{ij}(y_1,\ldots,y_{m-1},0) \mbox{ for } (i,j)
\mbox{ with }|x_i - x_j|=r \bigr\}.
\]
Indeed
\begin{eqnarray*}&& \max_{\{i<j: |x_i - x_j|=r\}}
g_{ij}(y) -\max_{\{i<j: |x_i - x_j|=r\}} g_{ij}
\bigl(y'\bigr)\\
&&\quad = g_{i_0j_0}(y) -\max
_{\{i<j: |x_{i} - x_{j}|=r\}} g_{ij}\bigl(y'\bigr) \qquad\mbox{for
some } i_0,j_0
\\
&&\quad\le g_{i_0j_0}(y) - g_{i_0j_0}\bigl(y'
\bigr) \le C_{i_0j_0}\bigl |y-y'\bigr|
\\
&&\quad\le \Bigl(\max_{\{i<j: |x_i - x_j|=r\}} C_{ij} \Bigr)
\bigl|y-y'\bigr|.
\end{eqnarray*}
Hence, $D$ is a Lipschitz domain.
\end{pf}
\end{appendix}

\section*{Acknowledgements}
The authors are very grateful to Giovanni Conforti, who realised the
simulations presented in the Figures \ref{fig1}--\ref{fig3}, see
Section~\ref{sec2}.

The first author was partially supported by the ANR project EVOL
(ANR-08-BLAN-0242).

The second and fourth authors gratefully acknowledge the support of the
Deutschen Akademischen Austauschdienst,
the BMBF and the french national agency CampusFrance through the
PROCOPE Project \emph{Structures g\'eom\'etriques
al\'eatoires en espace et en temps} (N\textsuperscript{o} 26491SJ) --
\emph{Strukturanalyse zuf\"alliger Geometrie in Raum und Zeit} (ID
54366242).

The third author acknowledges the Berlin Mathematical School and the
Research Training Group 1845
\emph{Stochastic Analysis with Applications in Biology, Finance and
Physics} for their partial support.

%

\printhistory

\begin{thebibliography}{39}

\bibitem{BCG}
%
\begin{barticle}[mr]
\bauthor{\bsnm{Bakry},~\bfnm{Dominique}\binits{D.}},
\bauthor{\bsnm{Cattiaux},~\bfnm{Patrick}\binits{P.}} \AND
\bauthor{\bsnm{Guillin},~\bfnm{Arnaud}\binits{A.}}
(\byear{2008}).
\btitle{Rate of convergence for ergodic continuous {M}arkov processes:
{L}yapunov versus {P}oincar\'e}.
\bjournal{J. Funct. Anal.}
\bvolume{254}
\bpages{727--759}.
\bid{doi={10.1016/j.jfa.2007.11.002}, issn={0022-1236}, mr={2381160}}
\end{barticle}
%
\bptok{imsref}%
\endbibitem

\bibitem{BassHsu}
%
\begin{barticle}[mr]
\bauthor{\bsnm{Bass},~\bfnm{Richard~F.}\binits{R.F.}} \AND
\bauthor{\bsnm{Hsu},~\bfnm{Pei}\binits{P.}}
(\byear{1990}).
\btitle{The semimartingale structure of reflecting {B}rownian motion}.
\bjournal{Proc. Amer. Math. Soc.}
\bvolume{108}
\bpages{1007--1010}.
\bid{doi={10.2307/2047960}, issn={0002-9939}, mr={1007487}}
\end{barticle}
%
\bptok{imsref}%
\endbibitem

\bibitem{BassHsu2}
%
\begin{barticle}[mr]
\bauthor{\bsnm{Bass},~\bfnm{Richard~F.}\binits{R.F.}} \AND
\bauthor{\bsnm{Hsu},~\bfnm{Pei}\binits{P.}}
(\byear{1991}).
\btitle{Some potential theory for reflecting {B}rownian motion in H\"
older and {L}ipschitz domains}.
\bjournal{Ann. Probab.}
\bvolume{19}
\bpages{486--508}.
\bid{issn={0091-1798}, mr={1106272}}
\end{barticle}
%
\bptok{imsref}%
\endbibitem

\bibitem{Bobsphere}
%
\begin{bincollection}[mr]
\bauthor{\bsnm{Bobkov},~\bfnm{S.~G.}\binits{S.G.}}
(\byear{2003}).
\btitle{Spectral gap and concentration for some spherically symmetric
probability measures}.
In \bbooktitle{Geometric Aspects of Functional Analysis}.
\bseries{Lecture Notes in Math.}
\bvolume{1807}
\bpages{37--43}.
\blocation{Berlin}:
\bpublisher{Springer}.
\bid{doi={10.1007/978-3-540-36428-3_4}, mr={2083386}}
\end{bincollection}
%
\bptok{imsref}%
\endbibitem

\bibitem{BCGM}
%
\begin{bmisc}[auto:parserefs-M02]
\bauthor{\bsnm{Boissard},~\bfnm{E.}\binits{E.}},
\bauthor{\bsnm{Cattiaux},~\bfnm{P.}\binits{P.}},
\bauthor{\bsnm{Guillin},~\bfnm{A.}\binits{A.}} \AND
\bauthor{\bsnm{Miclo},~\bfnm{L.}\binits{L.}}
(\byear{2013}).
\bhowpublished{Ornstein--Uhlenbeck pinball: I. Poincar\'e Inequalities
in a punctured domain.
Preprint. Available at \surl{http://perso.math.univ-toulouse.fr/cattiaux/}.}
\end{bmisc}
%
\bptok{imsref}%
\endbibitem

\bibitem{Boeroe}
%
\begin{bbook}[mr]
\bauthor{\bsnm{B{\"o}r{\"o}czky},~\bfnm{K{\'a}roly}\binits{K.}
\bsuffix{Jr.}}
(\byear{2004}).
\btitle{Finite Packing and Covering}.
\bseries{Cambridge Tracts in Mathematics}
\bvolume{154}.
\blocation{Cambridge}:
\bpublisher{Cambridge Univ. Press}.
\bid{doi={10.1017/CBO9780511546587}, mr={2078625}}
\end{bbook}
%
\bptok{imsref}%
\endbibitem

\bibitem{BoSa}
%
\begin{bbook}[mr]
\bauthor{\bsnm{Borodin},~\bfnm{Andrei~N.}\binits{A.N.}} \AND
\bauthor{\bsnm{Salminen},~\bfnm{Paavo}\binits{P.}}
(\byear{2002}).
\btitle{Handbook of {B}rownian Motion -- Facts and Formulae},
\bedition{2nd} ed.
\bseries{Probability and Its Applications}.
\blocation{Basel}:
\bpublisher{Birkh\"auser}.
\bid{doi={10.1007/978-3-0348-8163-0}, mr={1912205}}
\end{bbook}
%
\bptok{imsref}%
\endbibitem

\bibitem{Cx86}
%
\begin{barticle}[mr]
\bauthor{\bsnm{Cattiaux},~\bfnm{Patrick}\binits{P.}}
(\byear{1986}).
\btitle{Hypoellipticit\'e et hypoellipticit\'e partielle pour les
diffusions avec une condition fronti\`ere}.
\bjournal{Ann. Inst. Henri Poincar\'e Probab. Stat.}
\bvolume{22}
\bpages{67--112}.
\bid{issn={0246-0203}, mr={0838373}}
\end{barticle}
%
\bptok{imsref}%
\endbibitem

\bibitem{Cx87}
%
\begin{barticle}[mr]
\bauthor{\bsnm{Cattiaux},~\bfnm{Patrick}\binits{P.}}
(\byear{1987}).
\btitle{R\'egularit\'e au bord pour les densit\'es et les densit\'es
conditionnelles d'une diffusion r\'efl\'echie hypoelliptique}.
\bjournal{Stochastics}
\bvolume{20}
\bpages{309--340}.
\bid{doi={10.1080/17442508708833447}, issn={0090-9491}, mr={0885877}}
\end{barticle}
%
\bptok{imsref}%
\endbibitem

\bibitem{Cx92}
%
\begin{barticle}[mr]
\bauthor{\bsnm{Cattiaux},~\bfnm{Patrick}\binits{P.}}
(\byear{1992}).
\btitle{Stochastic calculus and degenerate boundary value problems}.
\bjournal{Ann. Inst. Fourier (Grenoble)}
\bvolume{42}
\bpages{541--624}.
\bid{issn={0373-0956}, mr={1182641}}
\end{barticle}
%
\bptok{imsref}%
\endbibitem

\bibitem{CGZ}
%
\begin{barticle}[mr]
\bauthor{\bsnm{Cattiaux},~\bfnm{Patrick}\binits{P.}},
\bauthor{\bsnm{Guillin},~\bfnm{Arnaud}\binits{A.}} \AND
\bauthor{\bsnm{Zitt},~\bfnm{Pierre~Andr{\'e}}\binits{P.A.}}
(\byear{2013}).
\btitle{Poincar\'e inequalities and hitting times}.
\bjournal{Ann. Inst. Henri Poincar\'e Probab. Stat.}
\bvolume{49}
\bpages{95--118}.
\bid{doi={10.1214/11-AIHP447}, issn={0246-0203}, mr={3060149}}
\end{barticle}
%
\bptok{imsref}%
\endbibitem

\bibitem{CFTYZ}
%
\begin{barticle}[mr]
\bauthor{\bsnm{Chen},~\bfnm{Z.-Q.}\binits{Z.-Q.}},
\bauthor{\bsnm{Fitzsimmons},~\bfnm{P.~J.}\binits{P.J.}},
\bauthor{\bsnm{Takeda},~\bfnm{M.}\binits{M.}},
\bauthor{\bsnm{Ying},~\bfnm{J.}\binits{J.}} \AND
\bauthor{\bsnm{Zhang},~\bfnm{T.-S.}\binits{T.-S.}}
(\byear{2004}).
\btitle{Absolute continuity of symmetric {M}arkov processes}.
\bjournal{Ann. Probab.}
\bvolume{32}
\bpages{2067--2098}.
\bid{doi={10.1214/009117904000000432}, issn={0091-1798}, mr={2073186}}
\end{barticle}
%
\bptok{imsref}%
\endbibitem

\bibitem{CFW}
%
\begin{barticle}[mr]
\bauthor{\bsnm{Chen},~\bfnm{Z.~Q.}\binits{Z.Q.}},
\bauthor{\bsnm{Fitzsimmons},~\bfnm{P.~J.}\binits{P.J.}} \AND
\bauthor{\bsnm{Williams},~\bfnm{R.~J.}\binits{R.J.}}
(\byear{1993}).
\btitle{Reflecting {B}rownian motions: Quasimartingales and strong
{C}accioppoli sets}.
\bjournal{Potential Anal.}
\bvolume{2}
\bpages{219--243}.
\bid{doi={10.1007/BF01048506}, issn={0926-2601}, mr={1245240}}
\end{barticle}
%
\bptok{imsref}%
\endbibitem

\bibitem{Chow}
%
\begin{barticle}[mr]
\bauthor{\bsnm{Chow},~\bfnm{Timothy~Y.}\binits{T.Y.}}
(\byear{1995}).
\btitle{Penny-packings with minimal second moments}.
\bjournal{Combinatorica}
\bvolume{15}
\bpages{151--158}.
\bid{doi={10.1007/BF01200751}, issn={0209-9683}, mr={1337349}}
\end{barticle}
%
\bptok{imsref}%
\endbibitem

\bibitem{ConwaySloane}
%
\begin{bbook}[mr]
\bauthor{\bsnm{Conway},~\bfnm{J.~H.}\binits{J.H.}} \AND
\bauthor{\bsnm{Sloane},~\bfnm{N.~J.~A.}\binits{N.J.A.}}
(\byear{1993}).
\btitle{Sphere Packings, Lattices and Groups},
\bedition{2nd} ed.
\bseries{Grundlehren der Mathematischen Wissenschaften}
\bvolume{290}.
\blocation{New York}:
\bpublisher{Springer}.
\bid{doi={10.1007/978-1-4757-2249-9}, mr={1194619}}
\end{bbook}
%
\bptok{imsref}%
\endbibitem

\bibitem{DLM}
%
\begin{barticle}[mr]
\bauthor{\bsnm{Diaconis},~\bfnm{Persi}\binits{P.}},
\bauthor{\bsnm{Lebeau},~\bfnm{Gilles}\binits{G.}} \AND
\bauthor{\bsnm{Michel},~\bfnm{Laurent}\binits{L.}}
(\byear{2011}).
\btitle{Geometric analysis for the Metropolis algorithm on {L}ipschitz
domains}.
\bjournal{Invent. Math.}
\bvolume{185}
\bpages{239--281}.
\bid{doi={10.1007/s00222-010-0303-6}, issn={0020-9910}, mr={2819161}}
\end{barticle}
%
\bptok{imsref}%
\endbibitem

\bibitem{FradonGlob}
%
\begin{barticle}[mr]
\bauthor{\bsnm{Fradon},~\bfnm{Myriam}\binits{M.}}
(\byear{2010}).
\btitle{Brownian dynamics of globules}.
\bjournal{Electron. J. Probab.}
\bvolume{15}
\bpages{142--161}.
\bid{doi={10.1214/EJP.v15-739}, issn={1083-6489}, mr={2594875}}
\end{barticle}
%
\bptok{imsref}%
\endbibitem

\bibitem{FR1}
%
\begin{barticle}[mr]
\bauthor{\bsnm{Fradon},~\bfnm{Myriam}\binits{M.}} \AND
\bauthor{\bsnm{R{\oe}lly},~\bfnm{Sylvie}\binits{S.}}
(\byear{2000}).
\btitle{Infinite-dimensional diffusion processes with singular interaction}.
\bjournal{Bull. Sci. Math.}
\bvolume{124}
\bpages{287--318}.
\bid{doi={10.1016/S0007-4497(00)00136-6}, issn={0007-4497}, mr={1771938}}
\end{barticle}
%
\bptok{imsref}%
\endbibitem

\bibitem{FR3}
%
\begin{barticle}[mr]
\bauthor{\bsnm{Fradon},~\bfnm{Myriam}\binits{M.}} \AND
\bauthor{\bsnm{R{\oe}lly},~\bfnm{Sylvie}\binits{S.}}
(\byear{2006}).
\btitle{Infinite system of {B}rownian balls: Equilibrium measures are
canonical {G}ibbs}.
\bjournal{Stoch. Dyn.}
\bvolume{6}
\bpages{97--122}.
\bid{doi={10.1142/S0219493706001669}, issn={0219-4937}, mr={2210683}}
\end{barticle}
%
\bptok{imsref}%
\endbibitem

\bibitem{FR2}
%
\begin{barticle}[mr]
\bauthor{\bsnm{Fradon},~\bfnm{Myriam}\binits{M.}} \AND
\bauthor{\bsnm{R{\oe}lly},~\bfnm{Sylvie}\binits{S.}}
(\byear{2007}).
\btitle{Infinite system of {B}rownian balls with interaction: The
non-reversible case}.
\bjournal{ESAIM Probab. Stat.}
\bvolume{11}
\bpages{55--79}.
\bid{doi={10.1051/ps:2007006}, issn={1292-8100}, mr={2299647}}
\end{barticle}
%
\bptok{imsref}%
\endbibitem

\bibitem{FradonRoellyGlob}
%
\begin{barticle}[mr]
\bauthor{\bsnm{Fradon},~\bfnm{Myriam}\binits{M.}} \AND
\bauthor{\bsnm{R{\oe}lly},~\bfnm{Sylvie}\binits{S.}}
(\byear{2010}).
\btitle{Infinitely many {B}rownian globules with {B}rownian radii}.
\bjournal{Stoch. Dyn.}
\bvolume{10}
\bpages{591--612}.
\bid{doi={10.1142/S021949371000311X}, issn={0219-4937}, mr={2740705}}
\end{barticle}
%
\bptok{imsref}%
\endbibitem

\bibitem{FRT}
%
\begin{barticle}[mr]
\bauthor{\bsnm{Fradon},~\bfnm{Myriam}\binits{M.}},
\bauthor{\bsnm{Roelly},~\bfnm{Sylvie}\binits{S.}} \AND
\bauthor{\bsnm{Tanemura},~\bfnm{Hideki}\binits{H.}}
(\byear{2000}).
\btitle{An infinite system of {B}rownian balls with infinite range
interaction}.
\bjournal{Stochastic Process. Appl.}
\bvolume{90}
\bpages{43--66}.
\bid{doi={10.1016/S0304-4149(00)00036-3}, issn={0304-4149}, mr={1787124}}
\end{barticle}
%
\bptok{imsref}%
\endbibitem

\bibitem{FOT}
%
\begin{bbook}[mr]
\bauthor{\bsnm{Fukushima},~\bfnm{Masatoshi}\binits{M.}},
\bauthor{\bsnm{{\=O}shima},~\bfnm{Y{\=o}ichi}\binits{Y.}} \AND
\bauthor{\bsnm{Takeda},~\bfnm{Masayoshi}\binits{M.}}
(\byear{1994}).
\btitle{Dirichlet Forms and Symmetric {M}arkov Processes}.
\bseries{De Gruyter Studies in Mathematics}
\bvolume{19}.
\blocation{Berlin}:
\bpublisher{de Gruyter}.
\bid{doi={10.1515/9783110889741}, mr={1303354}}
\end{bbook}
%
\bptok{imsref}%
\endbibitem

\bibitem{Fuktom}
%
\begin{barticle}[mr]
\bauthor{\bsnm{Fukushima},~\bfnm{Masatoshi}\binits{M.}} \AND
\bauthor{\bsnm{Tomisaki},~\bfnm{Matsuyo}\binits{M.}}
(\byear{1996}).
\btitle{Construction and decomposition of reflecting diffusions on
{L}ipschitz domains with H\"older cusps}.
\bjournal{Probab. Theory Related Fields}
\bvolume{106}
\bpages{521--557}.
\bid{doi={10.1007/s004400050074}, issn={0178-8051}, mr={1421991}}
\end{barticle}
%
\bptok{imsref}%
\endbibitem

\bibitem{Kinkladze82}
%
\begin{barticle}[mr]
\bauthor{\bsnm{Kinkladze},~\bfnm{G.~N.}\binits{G.N.}}
(\byear{1982}).
\btitle{A note on the structure of processes the measure of which is
absolutely continuous with respect to the {W}iener process modulus measure}.
\bjournal{Stochastics}
\bvolume{8}
\bpages{39--44}.
\bid{doi={10.1080/17442508208833226}, issn={0090-9491}, mr={0687044}}
\bptnote{check year}%
\end{barticle}
%
\bptok{imsref}%
\endbibitem

\bibitem{Kul2}
%
\begin{barticle}[mr]
\bauthor{\bsnm{Kulik},~\bfnm{Alexey~M.}\binits{A.M.}}
(\byear{2009}).
\btitle{Exponential ergodicity of the solutions to SDEs with a jump noise}.
\bjournal{Stochastic Process. Appl.}
\bvolume{119}
\bpages{602--632}.
\bid{doi={10.1016/j.spa.2008.02.006}, issn={0304-4149}, mr={2494006}}
\end{barticle}
%
\bptok{imsref}%
\endbibitem

\bibitem{Kul1}
%
\begin{barticle}[mr]
\bauthor{\bsnm{Kulik},~\bfnm{Alexey~M.}\binits{A.M.}}
(\byear{2011}).
\btitle{Asymptotic and spectral properties of exponentially {$\phi
$}-ergodic {M}arkov processes}.
\bjournal{Stochastic Process. Appl.}
\bvolume{121}
\bpages{1044--1075}.
\bid{doi={10.1016/j.spa.2011.01.007}, issn={0304-4149}, mr={2775106}}
\end{barticle}
%
\bptok{imsref}%
\endbibitem

\bibitem{Kul3}
%
\begin{barticle}[mr]
\bauthor{\bsnm{Kulik},~\bfnm{Alexey~M.}\binits{A.M.}}
(\byear{2011}).
\btitle{Poincar\'e inequality and exponential integrability of the
hitting times of a {M}arkov process}.
\bjournal{Theory Stoch. Process.}
\bvolume{17}
\bpages{71--80}.
\bid{issn={0321-3900}, mr={2934560}}
\end{barticle}
%
\bptok{imsref}%
\endbibitem

\bibitem{Linets}
%
\begin{barticle}[mr]
\bauthor{\bsnm{Linetsky},~\bfnm{Vadim}\binits{V.}}
(\byear{2005}).
\btitle{On the transition densities for reflected diffusions}.
\bjournal{Adv. in Appl. Probab.}
\bvolume{37}
\bpages{435--460}.
\bid{doi={10.1239/aap/1118858633}, issn={0001-8678}, mr={2144561}}
\end{barticle}
%
\bptok{imsref}%
\endbibitem

\bibitem{MT}
%
\begin{bbook}[mr]
\bauthor{\bsnm{Meyn},~\bfnm{S.~P.}\binits{S.P.}} \AND
\bauthor{\bsnm{Tweedie},~\bfnm{R.~L.}\binits{R.L.}}
(\byear{1993}).
\btitle{Markov Chains and Stochastic Stability}.
\bseries{Communications and Control Engineering Series}.
\blocation{London}:
\bpublisher{Springer}.
\bid{doi={10.1007/978-1-4471-3267-7}, mr={1287609}}
\end{bbook}
%
\bptok{imsref}%
\endbibitem

\bibitem{Osada96}
%
\begin{barticle}[mr]
\bauthor{\bsnm{Osada},~\bfnm{Hirofumi}\binits{H.}}
(\byear{1996}).
\btitle{Dirichlet form approach to infinite-dimensional {W}iener
processes with singular interactions}.
\bjournal{Comm. Math. Phys.}
\bvolume{176}
\bpages{117--131}.
\bid{issn={0010-3616}, mr={1372820}}
\end{barticle}
%
\bptok{imsref}%
\endbibitem

\bibitem{ST86}
%
\begin{barticle}[mr]
\bauthor{\bsnm{Saisho},~\bfnm{Yasumasa}\binits{Y.}} \AND
\bauthor{\bsnm{Tanaka},~\bfnm{Hiroshi}\binits{H.}}
(\byear{1986}).
\btitle{Stochastic differential equations for mutually reflecting
{B}rownian balls}.
\bjournal{Osaka J. Math.}
\bvolume{23}
\bpages{725--740}.
\bid{issn={0030-6126}, mr={0866273}}
\end{barticle}
%
\bptok{imsref}%
\endbibitem

\bibitem{ST87}
%
\begin{barticle}[mr]
\bauthor{\bsnm{Saisho},~\bfnm{Yasumasa}\binits{Y.}} \AND
\bauthor{\bsnm{Tanaka},~\bfnm{Hiroshi}\binits{H.}}
(\byear{1987}).
\btitle{On the symmetry of a reflecting {B}rownian motion defined by
{S}korohod's equation for a multidimensional domain}.
\bjournal{Tokyo J. Math.}
\bvolume{10}
\bpages{419--435}.
\bid{doi={10.3836/tjm/1270134524}, issn={0387-3870}, mr={0926253}}
\end{barticle}
%
\bptok{imsref}%
\endbibitem

\bibitem{SHDC}
%
\begin{barticle}[mr]
\bauthor{\bsnm{Sloane},~\bfnm{N.~J.~A.}\binits{N.J.A.}},
\bauthor{\bsnm{Hardin},~\bfnm{R.~H.}\binits{R.H.}},
\bauthor{\bsnm{Duff},~\bfnm{T.~D.~S.}\binits{T.D.S.}} \AND
\bauthor{\bsnm{Conway},~\bfnm{J.~H.}\binits{J.H.}}
(\byear{1995}).
\btitle{Minimal-energy clusters of hard spheres}.
\bjournal{Discrete Comput. Geom.}
\bvolume{14}
\bpages{237--259}.
\bid{doi={10.1007/BF02570704}, issn={0179-5376}, mr={1344734}}
\end{barticle}
%
\bptok{imsref}%
\endbibitem

\bibitem{Tanemura96}
%
\begin{barticle}[mr]
\bauthor{\bsnm{Tanemura},~\bfnm{Hideki}\binits{H.}}
(\byear{1996}).
\btitle{A system of infinitely many mutually reflecting {B}rownian
balls in {${\mathbf R}^d$}}.
\bjournal{Probab. Theory Related Fields}
\bvolume{104}
\bpages{399--426}.
\bid{doi={10.1007/BF01213687}, issn={0178-8051}, mr={1376344}}
\end{barticle}
%
\bptok{imsref}%
\endbibitem

\bibitem{Tanemura97}
%
\begin{barticle}[mr]
\bauthor{\bsnm{Tanemura},~\bfnm{Hideki}\binits{H.}}
(\byear{1997}).
\btitle{Uniqueness of {D}irichlet forms associated with systems of
infinitely many {B}rownian balls in {$\mathbf{R}^d$}}.
\bjournal{Probab. Theory Related Fields}
\bvolume{109}
\bpages{275--299}.
\bid{doi={10.1007/s004400050133}, issn={0178-8051}, mr={1477652}}
\end{barticle}
%
\bptok{imsref}%
\endbibitem

\bibitem{Tem}
%
\begin{barticle}[mr]
\bauthor{\bsnm{Temesv{\'a}ri},~\bfnm{{\'A}gota}\binits{{\'A}.}}
(\byear{1974}).
\btitle{On the extremum of power sums of distances}.
\bjournal{Mat. Lapok (N.S.)}
\bvolume{25}
\bpages{329--342}.
\bid{issn={0025-519X}, mr={0514018}}
\end{barticle}
%
\bptok{imsref}%
\endbibitem

\bibitem{Veret}
%
\begin{barticle}[mr]
\bauthor{\bsnm{Veretennikov},~\bfnm{A.~Y.}\binits{A.Y.}}
(\byear{1987}).
\btitle{Estimates of the mixing rate for stochastic equations}.
\bjournal{Teor. Veroyatn. Primen.}
\bvolume{32}
\bpages{299--308}.
\bid{issn={0040-361X}, mr={0902757}}
\end{barticle}
%
\bptok{imsref}%
\endbibitem

\bibitem{WG}
%
\begin{barticle}[mr]
\bauthor{\bsnm{Ward},~\bfnm{Amy~R.}\binits{A.R.}} \AND
\bauthor{\bsnm{Glynn},~\bfnm{Peter~W.}\binits{P.W.}}
(\byear{2003}).
\btitle{Properties of the reflected {O}rnstein--{U}hlenbeck process}.
\bjournal{Queueing Syst.}
\bvolume{44}
\bpages{109--123}.
\bid{doi={10.1023/A:1024403704190}, issn={0257-0130}, mr={1993278}}
\end{barticle}
%
\bptok{imsref}%
\endbibitem
\end{thebibliography}
\end{document}